\documentclass[journal,onecolumn,draftclsnofoot,11pt]{IEEEtran}

\usepackage{amsthm}
\usepackage{amsbsy}
\usepackage{amsmath}
\usepackage{amssymb}
\usepackage{latexsym}
\usepackage{ifthen}
\usepackage{subfigure}
\usepackage{turnstile}
\usepackage{mathtools}
\usepackage{booktabs}
\usepackage{balance}
\usepackage{paralist}
\usepackage{cite}
\usepackage{color}
\usepackage{url}
\usepackage{caption}

\newtheorem{Theorem}{Theorem}[section]
\newtheorem{Lemma}[Theorem]{Lemma}
\newtheorem{proposition}[Theorem]{Proposition}

\newtheorem{corollary}[Theorem]{Corollary}
\newtheorem{Remark}[Theorem]{Remark}
\theoremstyle{definition}
\newtheorem{assump}{Assumption}

\newenvironment{myassump}[2][]
  {\renewcommand\theassump{#2}\begin{assump}[#1]}
  {\end{assump}}

\def \cisprt{\mathcal{CI}SPRT}

\mathtoolsset{showonlyrefs}

\begin{document}

\title{Distributed Sequential Detection for Gaussian Shift-in-Mean Hypothesis Testing}
\author{Anit~Kumar~Sahu,~\IEEEmembership{Student Member,~IEEE} and Soummya~Kar,~\IEEEmembership{Member,~IEEE}
\thanks{The authors are with the Department of Electrical and Computer Engineering, Carnegie Mellon University, Pittsburgh, PA 15213, USA (email: anits@andrew.cmu.edu, soummyak@andrew.cmu.edu). This work was partially supported by NSF Grant ECCS-1306128.}}
\maketitle
\begin{abstract}
This paper studies the problem of sequential Gaussian shift-in-mean hypothesis testing in a distributed multi-agent network. A sequential probability ratio test (SPRT) type algorithm in a distributed framework of the \emph{consensus}+\emph{innovations} form is proposed, in which the agents update their decision statistics by simultaneously processing latest observations (innovations) sensed sequentially over time and information obtained from neighboring agents (consensus).  For each pre-specified set of type I and type II error probabilities, local decision parameters are derived which ensure that the algorithm achieves the desired error performance and terminates in finite time almost surely (a.s.) at each network agent. Large deviation exponents for the tail probabilities of the agent stopping time distributions are obtained and it is shown that asymptotically (in the number of agents or in the high signal-to-noise-ratio regime) these exponents associated with the distributed algorithm approach that of the optimal centralized detector. The expected stopping time for the proposed algorithm at each network agent is evaluated and is benchmarked with respect to the optimal centralized algorithm. The efficiency of the proposed algorithm in the sense of the expected stopping times is characterized in terms of network connectivity. Finally, simulation studies are presented which illustrate and verify the analytical findings.
\end{abstract}

\begin{IEEEkeywords}
Distributed Detection, Multi-agent Networks, Consensus, Sequential Probability Ratio Tests, Large Deviations
\end{IEEEkeywords}

\section{Introduction}
\subsection{Background and Motivation}
The focus of this paper is on sequential simple hypothesis testing in multi-agent networks in which the goal is to detect the (binary) state of the environment based on observations at the agents. By sequential we mean, instead of considering fixed sample size hypothesis tests in which the objective is to minimize the probabilities of decision error (the false alarm and the miss) based on a given deterministic number of samples or observation data collected by the network agents, we are interested in the design of testing procedures that in the \emph{quickest} time or using the \emph{minimal} amount of sensed data samples at the agents can distinguish between the two hypotheses with guaranteed accuracy given in terms of pre-specified tolerances on false alarm and miss probabilities. The motivation behind studying sequential as opposed to fixed sample size testing is that in most practical agent networking scenarios, especially in applications that are time-sensitive and/or resource constrained, the priority is to achieve inference as quickly as possible by expending the minimal amount of resources (data samples, sensing energy and communication). Furthermore, we focus on distributed application environments which are devoid of fusion centers\footnote{By fusion center or center, we mean a hypothetical decision-making architecture in which a (central) entity has access to all agent observations at all times and/or is responsible for decision-making on behalf of the agents.} and in which inter-agent collaboration or information exchange is limited to a pre-assigned, possibly sparse, communication structure.

Under rather generic assumptions on the agent observation models, it is well-known that in a (hypothetical) centralized scenario or one in which inter-agent communication is all-to-all corresponding to a complete communication graph, the sequential probability ratio test (SPRT) (\hspace{-0.5pt}\cite{wald1945sequential}) turns out to be the optimal procedure for sequential testing of binary hypotheses; specifically, the SPRT minimizes the expected detection time (and hence the number of agent observation samples that need to be processed) while achieving requisite error performance in terms of specified probability of false alarm ($\alpha$) and probability of miss ($\beta$) tolerances. The SPRT and its variants have been applied in various contexts, see, for example, spectrum sensing in cognitive radio networks (\hspace{-0.5pt}\cite{choi2009sequential,jayaprakasam2009sequential,chaudhari2009autocorrelation}), target tracking~\cite{blostein1994sequential}, to name a few. However, the SPRT, in the current multi-agent context, would require computing a (centralized) decision statistic at all times, which, in turn, would either require all-to-all communication among the agents or access to the entire network data at all times at a fusion center. In contrast, restricted by a pre-assigned possibly sparse collaboration structure among the agents, in this paper we present and characterize a distributed sequential detection algorithm, the $\cisprt$, based on the \emph{consensus}+\emph{innovations} approach (see, for example \cite{kar-asilomar-linest-2008,KarMouraRamanan-Est-2008}). Specifically, focusing on a setting in which the agent observations over time are conditionally Gaussian and independent and identically distributed (i.i.d.), we study the $\cisprt$ sequential detection procedure in which each network agent maintains a local (scalar) test statistic which is updated over time by simultaneously assimilating the test statistics of neighboring agents at the previous time instant (a consensus potential) and the most recent observations (innovations) obtained by the agent and its neighbors. Also, similar in spirit to the (centralized) SPRT, each agent chooses two (local) threshold parameters (design choices) and the test termination at an agent (and subsequent agent decision on the hypotheses) is determined by whether the local test statistic at the agent lies in the interval defined by the thresholds or not. This justifies the nomenclature that the $\cisprt$ is a distributed SPRT type algorithm of the consensus+innovations form.
The main contributions of this paper are as follows:

\noindent\textbf{Main Contribution 1: Finite Stopping Property.} We show that, given any value of probability of false alarm $\alpha$ and probability of miss $\beta$, the $\cisprt$ algorithm can be designed such that each agent achieves the specified error performance metrics and the test procedure terminates in finite time almost surely (a.s.) at each agent. We derive closed form expressions for the local threshold parameters at the agents as functions of $\alpha$ and $\beta$ which ensures that the $\cisprt$ achieves the above property.\\
\noindent\textbf{Main Contribution 2: Asymptotic Characterization.} By characterizing the stopping time distribution of the $\cisprt$ at each network agent, we compute large deviations decay exponents of the stopping time tail probabilities at each agent, and show that the large deviations exponent of the $\cisprt$ approaches that of the optimal centralized in the asymptotics of $N$, where $N$ denotes the number of agents in the network. In the asymptotics of vanishing error metrics (i.e., as $\alpha,\beta\rightarrow 0$), we quantify the ratio of the expected stopping time $T_{d,i}(\alpha,\beta)$ for reaching a decision at an agent $i$ through the $\cisprt$ algorithm and the expected stopping time $T_{c}(\alpha,\beta)$ for reaching a decision by the optimal centralized (SPRT) procedure, i.e., the quantity $\frac{\mathbb{E}[T_{d,i}(\alpha,\beta)]}{\mathbb{E}[T_{c}(\alpha,\beta)]}$, which in turn is a metric of efficiency of the proposed algorithm as a function of the network connectivity. In particular, we show that the efficiency of the proposed $\cisprt$ algorithm in terms of the ratio $\frac{\mathbb{E}[T_{d,i}(\alpha,\beta)]}{\mathbb{E}[T_{c}(\alpha,\beta)]}$ is upper bounded by a constant which is a function of the network connectivity and can be made close to one by choosing the network connectivity appropriately, thus establishing the benefits of inter-agent collaboration in the current context.\\
\textbf{Related Work.}
Detection schemes in multi-agent networks which involve fusion centers, where all agents in the network transmit their local measurements, local decisions or local likelihood ratios to a fusion agent which subsequently makes the final decision (see, for example, \cite{chamberland2003decentralized, tsitsiklis1993decentralized, blum1997distributed, veeravalli1993decentralized}) have been well studied.
Consensus-based approaches for fully distributed but single snapshot processing, i.e., in which the agents first collect their observations possibly over a long time horizon and then deploy a consensus-type protocol~\cite{jadbabailinmorse03,olfatisaberfaxmurray07,dimakis2010gossip} to obtain distributed information fusion and decision-making have also been explored, see, for instance,~\cite{kar2008topology, kar2007consensus}. Generalizations and variants of this framework have been developed, see for instance~\cite{zhang2014asymptotically} which proposes truncated versions of optimal testing procedures to facilitate efficient distributed computation using consensus; scenarios involving distributed information processing where some of the agents might be faulty or there is imperfect model information (see, for example, \cite{zhou2011robust, zhou2012distributed}) have also been studied. More relevant to the current context are distributed detection techniques that like the $\cisprt$ procedure perform simultaneous assimilation of neighborhood decision-statistics and local agent observations in the same time step, see, in particular, the running consensus approach~\cite{braca2008enforcing,braca2010asymptotic}, the diffusion approach~\cite{cattivelli2009diffusion,cattivelli2009distributed,cattivelli2011distributed} and the consensus+innovations approach~\cite{bajovic2011distributed ,jakovetic2012distributed,kar2011distributed}. These works address important questions in fixed (but possibly large) sample size distributed hypothesis testing, including asymptotic characterization of detection errors~\cite{braca2010asymptotic,cattivelli2011distributed}, fundamental performance limits as characterized by large deviations decay of detection error probabilities in generic nonlinear observation models and random networks~\cite{bajovic2011distributed ,jakovetic2012distributed}, and detection with noisy communication links~\cite{kar2011distributed}. A continuous time version of the running consensus approach~\cite{braca2010asymptotic} was studied in~\cite{srivastava2014collective} recently with implications on sequential distributed detection; specifically, asymptotic properties of the continuous time decision statistics were obtained and in the regime of large thresholds bounds on expected decision time and error probability rates were derived. However, there is a fundamental difference between the mostly fixed or large sample size procedures discussed above and the proposed $\cisprt$ sequential detection procedure -- technically speaking, the former focuses on analyzing the probability distributions of the detection errors as a function of the sample size and/or specified thresholds, whereas, in this paper, we design thresholds, stopping criteria and characterize the probability distributions of the (random) stopping times of sequential distributed procedures that aim to achieve quickest detection given specified tolerances on the detection errors. Addressing the latter requires novel technical machinery in the design and analysis of dynamic distributed inference procedures which we develop in this paper.

We also contrast our work with sequential detection approaches based in other types of multi-agent networking scenarios. In the context of decentralized sequential testing in multi-agent networks, fundamental methodological advances have been reported, see, for instance,~\cite{hashemi1989decentralized,mei2008asymptotic,nayyar2011sequential,veeravalli1994decentralized,fellouris2011decentralized,poor2009quickest,tartakovsky2004change,tartakovsky2000sequential}, which address very general models and setups. These works involve fusion center based processing where all agents in the network either transmit their local decisions, measurements or their quantized versions to a fusion center. In contrast, in this paper we restrict attention to Gaussian binary testing models only, but focus on a fully distributed paradigm in which there is no fusion center and inter-agent collaboration is limited to a pre-assigned, possibly sparse, agent-to-agent local interaction graph.\\

\textbf{Paper Organization :}
We briefly summarize the organization of the rest of the paper. Section~\ref{subsec:not} presents notation to be used throughout the paper. The sensing models and the abstract problem formulation are stated and discussed in Section~\ref{subsec:sys_model}. Section~\ref{subsec:prel_seq_test} presents preliminaries on centralized sequential detection and motivates the distributed setup pursued in this paper. Section~\ref{dsd} presents the $\cisprt$ algorithm. The main results of the paper are stated in Section~\ref{sec:main_res} which includes the derivation of the thresholds for the $\cisprt$ algorithm, the stopping time distribution for the $\cisprt$ algorithm and the key technical ingredients concerning the asymptotic properties and large deviation analysis for the stopping time distributions of the centralized and distributed setups. It also includes the characterization of the expected stopping times of the $\cisprt$ algorithm and its centralized counterpart in asymptotics of vanishing error metrics. Section~\ref{sec:sim} presents simulation studies. The proofs of the main results appear in Section~\ref{sec:proof_res}, whereas, Section~\ref{sec:conc} concludes the paper.

\subsection{Notation}
\label{subsec:not}
We denote by~$\mathbb{R}$ the set of reals, $\mathbb{R}_{+}$ the set of non-negative reals, and by~$\mathbb{R}^{k}$ the $k$-dimensional Euclidean space.
The set of $k\times k$ real matrices is denoted by $\mathbb{R}^{k\times k}$. The set of integers is denoted by $\mathbb{Z}$, whereas, $\mathbb{Z}_{+}$ denotes the subset of non-negative integers and $\overline{\mathbb{Z}}_{+}=\mathbb{Z}_{+}\cup\{\infty\}$. We denote vectors and matrices by bold faced characters. We denote by $A_{ij}$ or $[A]_{ij}$ the $(i,j)$th entry of a matrix $\mathbf{A}$; $a_{i}$ or $[a]_{i}$ the $i$th entry of a vector $\mathbf{a}$. The symbols $\mathbf{I}$ and $\mathbf{0}$ are used to denote the $k\times k$ identity matrix and the $k\times p$ zero matrix respectively, the dimensions being clear from the context. We denote by $\mathbf{e}_{i}$ the $i$th column of $\mathbf{I}$. The symbol $\top$ denotes matrix transpose. The $k\times k$ matrix $\mathbf{J}=\frac{1}{k}\mathbf{1}\mathbf{1^{\top}}$ where $\mathbf{1}$ denotes the $k\times 1$ vector of ones. The operator $\|\cdot\|$ applied to a vector denotes the standard Euclidean $\mathcal{L}_{2}$ norm, while applied to matrices it denotes the induced $\mathcal{L}_{2}$ norm, which is equivalent to the spectral radius for symmetric matrices. All the logarithms in the paper are with respect to base $e$ and represented as $\log(\cdot)$. Expectation is denoted by $\mathbb{E}[\cdot]$ and $\mathbb{E}_{\theta}[\cdot]$ denotes expectation conditioned on hypothesis $H_{\theta}$ for $\theta \in \{0,1\}$. $\mathbb{P}(\cdot)$ denotes the probability of an event and $\mathbb{P}_{\theta}( . )$ denotes the probability of the event conditioned on hypothesis $H_{\theta}$ for $\theta \in \{0,1\}$. $\mathbb{Q}( . )$ denotes the $\mathbb{Q}$-function which calculates the right tail probability of a normal distribution and is given by $\mathbb{Q}(x)=\frac{1}{\sqrt{2\pi}}\int_{x}^{\infty}e^{-\frac{u^{2}}{2}}du$, $x \in \mathbb{R}$. We will use the following property of the $\mathbb{Q}(.)$ function, namely for any $x >0$, $\mathbb{Q}(x) \le \frac{1}{2}e^{-\frac{x^{2}}{2}}$. For deterministic $\mathbb{R}_{+}$-valued sequences $\{a_{t}\}$ and $\{b_{t}\}$, the notation $a_{t}=O(b_{t})$ denotes the existence of a constant $c>0$ such that $a_{t}\leq cb_{t}$ for all $t$ sufficiently large, whereas, $a_{t}=o(b_{t})$ indicates $a_{t}/b_{t}\rightarrow 0$ as $t\rightarrow\infty$.\\

\noindent{\bf Spectral Graph Theory.} For an undirected graph $G=(V, E)$, $V$ denotes the set of agents or vertices with cardinality $|V|=N$, and $E$ the set of edges with $|E|=M$. The unordered pair $(i,j) \in E$ if there exists an edge between agents $i$ and $j$. We only consider simple graphs, i.e. graphs devoid of self loops and multiple edges. A path between agents $i$ and $j$ of length $m$ is a sequence ($i=p_{0},p_{1},\hdots,p_{m}=j)$ of vertices, such that $(p_{n}, p_{n+1})\in E$, $0\le n \le m-1$. A graph is connected if there exists a path between all the possible agent pairings.
The neighborhood of an agent $i$ is given by $\Omega_{i}=\{j \in V~|~(i,j) \in E\}$. The degree of agent $i$ is given by the cardinality $d_{i}=|\Omega_{i}|$. The structure of the graph may be equivalently represented by the symmetric $N\times N$ adjacency matrix $\mathbf{A}=[A_{ij}]$, where $A_{ij}=1$ if $(i,j) \in E$, and $0$ otherwise. The degree matrix is represented by the diagonal matrix $\mathbf{D}=diag(d_{1}\hdots d_{N})$. The graph Laplacian matrix is represented by
\begin{align}
\label{eq:0.2}
\mathbf{L}=\mathbf{D}-\mathbf{A}.
\end{align}
The Laplacian is a positive semidefinite matrix, hence its eigenvalues can be sorted and represented in the following manner
\begin{align}
\label{eq:0.3}
0=\lambda_{1}(\mathbf{L})\le\lambda_{2}(\mathbf{L})\le \hdots \lambda_{N}(\mathbf{L}).
\end{align}
Furthermore, a graph is connected if and only if $\lambda_{2}(\mathbf{L})>0$ (see~\cite{chung1997spectral}~for instance). We stress that under our notation, $\lambda_{2}$, also known as the Fiedler value (see~\cite{chung1997spectral}), plays an important role because it acts as an indicator of whether the graph is connected or not.

\section{Problem Formulation}
\label{sec:prob_form}

\subsection{System Model}
\label{subsec:sys_model}
\noindent The $N$ agents deployed in the network decide on either of the two hypothesis $H_{0}$ and $H_{1}$. Each agent $i$ at (discrete) time $t$ makes a scalar observation $y_{i}(t)$ of the form
\begin{align}
\label{eq:0.6}
Under~~ H_{\theta} : y_{i}(t)= \mu_{\theta}+n_{i}(t),~~\theta=0,1.
\end{align}

\noindent For the rest of the paper we consider $\mu_{1}=\mu$ and $\mu_{0}=-\mu$, and assume that, the agent observation noise processes are independent and identically distributed (i.i.d.) Gaussian processes under both hypotheses formalized as follows:

\begin{myassump}{A1}
\label{as:1}
\emph{For each agent $i$ the noise sequence $\{n_{i}(t)\}$ is i.i.d. Gaussian with mean zero and variance $\sigma^{2}$ under both $H_{0}$ and $H_{1}$. The noise sequences are also spatially uncorrelated, i.e., $\mathbb{E}_{\theta}[n_{i}(t)n_{j}(t)]=\mathbf{0}$ for all $i\neq j$ and $\theta\in\{0,1\}$.}
\end{myassump}

\noindent Collect the $y_{i}(t)$'s, $i=1,2,\hdots N$ into the $N \times 1$ vector $\mathbf{y}(t)=(y_{1}(t),\hdots,y_{N}(t))^{\top}$ and the  $n_{i}(t)$'s, $i=1,2,\hdots N$ into the $N \times 1$ vector $\mathbf{n}(t)=(n_{1}(t),\hdots,n_{N}(t))^{\top}$.

\noindent The log-likelihood ratio at the $i$-th sensor at time index $t$ is calculated as follows:-
\begin{align}
\label{eq:-0.5}
\eta_{i}(t)=\frac{f_{1}(y_{i}(t))}{f_{0}(y_{i}(t))}=\frac{2\mu y_{i}(t)}{\sigma^{2}},
\end{align}
where $f_{0}(\cdot)$ and $f_{1}(\cdot)$ denote the probability distribution functions (p.d.f.s) of $y_{i}(t)$ under $H_{0}$ and $H_{1}$ respectively.

\noindent We note that,
\begin{align}
\label{eq:f8}
\eta_{i}(t)\sim
\begin{cases}
\mathcal{N}(m,2m),& H=H_{1}\\
\mathcal{N}(-m,2m),& H=H_{0},
\end{cases}
\end{align}
where $\mathcal{N}(\cdot)$ denotes the Gaussian p.d.f. and $m=\frac{2\mu^{2}}{\sigma^{2}}$.
The Kullback-Leibler divergence at each agent is given by
\begin{align}
\label{eq:snr}
KL=m.
\end{align}

\subsection{Sequential Hypothesis Testing -- Centralized or All-To-All Communication Scenario}
\label{subsec:prel_seq_test} We start by reviewing concepts and results from (centralized) sequential hypothesis testing theory, see~\cite{chernoff1972sequential} for example, to motivate our distributed hypothesis testing setup. Broadly speaking, the goal of sequential simple hypothesis testing is as follows: given pre-specified constraints on the error metrics, i.e., upper bounds $\alpha$ and $\beta$ on the probability of false alarm $\mathbb{P}_{FA}$ and probability of miss $\mathbb{P}_{M}$, the decision-maker keeps on collecting observations sequentially over time to decide on the hypotheses $H_{1}$ or $H_{0}$, i.e., which one is true; the decision-maker also has a \emph{stopping criterion} or \emph{stopping rule} based on which it decides at each time (sampling) instant whether to continue sampling or terminate the testing procedure. Finally, after termination, a (binary) decision is computed as to which hypothesis is in force based on all the obtained data. A sequential testing procedure is said to be \emph{admissible} if the stopping criterion, i.e., the decision whether to continue observation collection or not, at each instant is determined solely on the basis of observations collected thus far. Naturally, from a resource optimization viewpoint, the decision-maker seeks to design the sequential procedure (or equivalently the stopping criterion) that minimizes the expected number of observation samples (or equivalently time) required to achieve a decision with probabilities of false alarm and miss upper bounded by $\alpha$ and $\beta$ respectively. To formalize in the current context, first consider a setup in which inter-agent communication is all-to-all, i.e., at each time instant each agent has access to all the sensed data of all other agents. In this complete network scenario, each agent behaves like a (hypothetical) center and the information available at any agent $n$ at time $t$ is the sum-total of network observations till $t$, formalized by the $\sigma$-algebra~\cite{jacod1987limit}
\begin{align}
\label{eq:0.01a}
\mathcal{G}_{c}(t)=\sigma\left\{y_{i}(s),~\forall i=1,2,\hdots N~\mbox{and}~\forall 1\le s \le t\right\}.
\end{align}
An admissible test $D_{c}$ consists of a stopping criteria, where at each time $t$ the agents' (or the center in this case) decision to stop or continue taking observations is adapted to (or measurable with respect to) the $\sigma$-algebra $\mathcal{G}_{c}(t)$. Denote by $T_{D_{c}}$ the termination time of $D_{c}$, a random time taking values in $\mathbb{Z}_{+}\cup\{\infty\}$. Formally, by the above notion of admissibility, the random time $T_{D_{c}}$ is necessarily a stopping time with respect to (w.r.t.) the filtration $\{\mathcal{G}_{c}(t)\}$, see~\cite{jacod1987limit}, and, in this paper, we restrict attention to tests $D_{c}$ that terminate in finite time a.s., i.e., $T_{D_{c}}$ takes values in $\mathbb{Z}_{+}$ a.s. Denote by $\mathbb{E}_{\theta}[T_{D_{c}}]$ the expectation of $T_{D_{c}}$ under $H_{\theta}$, $\theta=0,1$, and $\widehat{H}_{D_{c}}\in\{0,1\}$ the decision obtained after termination. (Note that, assuming $T_{D_{c}}$ is finite a.s., the random variable $\widehat{H}_{D_{c}}$ is measurable w.r.t. the stopped $\sigma$-algebra $\mathcal{G}_{T_{D_{c}}}$.) Let $\mathbb{P}_{FA}^{D_{c}}$ and $\mathbb{P}_{M}^{D_{c}}$ denote the associated probabilities of false alarm and miss respectively, i.e.,
\begin{equation}
\label{admiss:FAM} \mathbb{P}_{FA}^{D_{c}}=\mathbb{P}_{0}\left(\widehat{H}_{D_{c}}=1\right)~~\mbox{and}~~\mathbb{P}_{M}^{D_{c}}=\mathbb{P}_{1}\left(\widehat{H}_{D_{c}}=0\right).
\end{equation}
Now, denoting by $\mathcal{D}_{c}$ the class of all such (centralized) admissible tests, the goal in sequential hypothesis testing is to obtain a test in $\mathcal{D}_{c}$ that minimizes the expected stopping time subject to attaining specified error constraints. Formally, we aim to solve\footnote{Note, in~\eqref{eq:0.01} the objective is to minimize the expected stopping time under hypothesis $H_{1}$. Alternatively, we might be interested in minimizing $\mathbb{E}_{0}[T_{D_{c}}]$ over all admissible tests; similarly, in a Bayesian setup with prior probabilities $p_{0}$ and $p_{1}$ on $H_{0}$ and $H_{1}$ respectively, the objective would consist of minimizing the overall expected stopping time $p_{0}\mathbb{E}_{0}[T_{D_{c}}]+p_{1}\mathbb{E}_{1}[T_{D_{c}}]$. However, it turns out that, in the current context, the Wald's SPRT~\cite{wald1948optimum} (to be discussed soon) can be designed to minimize each of the above criteria. Hence, without loss of generality, we adopt $\mathbb{E}_{1}[T_{D_{c}}]$ as our test design objective and use it as a metric to determine the relative performance of tests.}
\begin{align}
\label{eq:0.01}
&\min\limits_{ D_{c}\in \mathcal{D}_{c}}\mathbb{E}_{1}[T_{D_{c}}], \nonumber\\
& \mbox{s.t.}~\mathbb{P}_{FA}^{D_{c}}\le\alpha, \mathbb{P}_{M}^{D_{c}}\le\beta,
\end{align}
for specified $\alpha$ and $\beta$. Before proceeding further, we make the following assumption:
\begin{myassump}{A2}
\label{as:2}
\emph{The pre-specified error metrics, i.e., $\alpha$ and $\beta$, satisfy $\alpha,\beta\in (0,1/2)$.}
\end{myassump}

Noting that the (centralized) Kullback-Leibler divergence, i.e., the divergence between the probability distributions induced on the joint observation space $\mathbf{y}(t)$ by the hypotheses $H_{1}$ and $H_{0}$, is $Nm$ where $m$ is defined in~\eqref{eq:f8}, we obtain (see~\cite{wald1948optimum}) for each $D_{c}\in\mathcal{D}_{c}$ that attains $\mathbb{P}_{FA}^{D_{c}}\le\alpha$ and $\mathbb{P}_{M}^{D_{c}}\le\beta$,
\begin{align}
\label{eq:edc300}
\mathbb{E}_{1}[T_{D_{c}}]\ge \mathcal{M}(\alpha,\beta),
\end{align}
where the universal lower bound $\mathcal{M}(\alpha,\beta)$ is given by
\begin{align}
\label{eq:edc3100}
\mathcal{M}(\alpha,\beta)=\frac{(1-\beta)\log(\frac{1-\beta}{\alpha})+\beta\log(\frac{\beta}{1-\alpha})}{Nm}.
\end{align}

\noindent{\bf Optimal (centralized) tests: Wald's SPRT.} We briefly review Wald's sequential probability ratio test (SPRT), see~\cite{wald1948optimum}, that is known to achieve optimality in~\eqref{eq:0.01}. To this end, denote by $S_{c}(t)$ (the centralized) test statistic
\begin{align}
\label{eq:cent_stat}
S_{c}(t)=\sum_{s=1}^{t}\frac{\mathbf{1}^{\top}}{N}\mathbf{\eta}(s),
\end{align}
where $\mathbf{\eta}(s)$ denotes the vector of log-likelihood ratios $\eta_{i}(s)$'s at the agents \footnote{Both the sum and average (over $N$) can be taken as the test statistics for the centralized detector. We divide by $N$ for notational simplicity, so that the centralized decision statistic update becomes a special case of the $\mathcal{CI}SPRT$ decision statistic update studied in Section~\ref{dsd}.}. The SPRT consists of a pair of thresholds (design parameters) $\gamma_{c}^{l}$ and $\gamma_{c}^{h}$, such that, at each time $t$, the decision to continue or terminate is determined on the basis of whether $S_{c}(t) \in [\gamma_{c}^{l},\gamma_{c}^{h}]$ or not. Formally, the stopping time of the SPRT is defined as follows:
\begin{align}
\label{eq:f5b}
 T_{c}=\inf \{t~|~S_{c}(t) \notin [\gamma_{c}^{l},\gamma_{c}^{h}]\}.
\end{align}
At $T_{c}$ the following decision rule is followed:
\begin{align}
\label{eq:f5c}
H=
\begin{cases}
H_{0}, & S_{c}(T_{c})\le \gamma_{c}^{l}\\
H_{1}, & S_{c}(T_{c})\ge \gamma_{c}^{h}.
\end{cases}
\end{align}
The optimality of the SPRT w.r.t. the formulation~\eqref{eq:0.01} is well-studied; in particular, in~\cite{wald1948optimum} it was shown that, for any specified $\alpha$ and $\beta$, there exist choices of thresholds $(\gamma_{c}^{l},\gamma_{c}^{h})$ such that the SPRT~\eqref{eq:f5b}-\eqref{eq:f5c} achieves the minimum in~\eqref{eq:0.01} among all possible admissible tests $D_{c}$ in $\mathcal{D}_{c}$.

For given $\alpha$ and $\beta$, exact analytical expressions of the optimal thresholds are intractable in general. A commonly used choice of thresholds, see~\cite{wald1945sequential}, is given by
\begin{align}
\label{eq:c1}
&\gamma_{c}^{h} = \log\big(\frac{1-\beta}{\alpha}\big)\nonumber\\
&\gamma_{c}^{l}= \log\big(\frac{\beta}{1-\alpha}\big),
\end{align}
which, although not strictly optimal in general, ensures that $\mathbb{P}_{FA}^{c}\le \alpha$ and $\mathbb{P}_{M}^{c}\le \beta$. (For SPRT procedures we denote by $\mathbb{P}_{FA}^{c}$ and $\mathbb{P}_{M}^{c}$ the associated probabilities of false alarm and miss respectively, which depend on the choice of thresholds used.) Nonetheless, the above choice~\eqref{eq:c1} yields \emph{close} to optimal behavior, and is, in fact, \emph{asymptotically} optimal; formally, supposing that $\alpha=\beta=\epsilon$, the SPRT with thresholds given by~\eqref{eq:c1} guarantees that (see~\cite{chernoff1972sequential})
\begin{equation}
\label{SPRT:cent:opt} \lim_{\epsilon\rightarrow 0}\frac{\mathbb{E}_{1}[T_{c}]}{\mathcal{M}(\epsilon,\epsilon)}=1,
\end{equation}
where $\mathcal{M}(\cdot)$ is defined in~\eqref{eq:edc3100}.
In the sequel, given a testing procedure $D_{c}\in\mathcal{D}_{c}$ and assuming $\alpha=\beta=\epsilon$, we will study the quantity $\limsup_{\epsilon\rightarrow 0}\left(\mathbb{E}_{1}[T_{D_{c}}]/\mathcal{M}(\epsilon,\epsilon)\right)$ as a measure of its efficiency. Also, by abusing notation, when $\alpha=\beta=\epsilon$, we will denote $\mathcal{M}(\epsilon)\doteq\mathcal{M}(\epsilon,\epsilon)$.

\subsection{Subclass of Distributed Tests}
\label{subsec:dist_tests} The SPRT~\eqref{eq:f5b}-\eqref{eq:f5c} requires computation of the statistic $S_{c}(t)$ (see~\eqref{eq:cent_stat}) at all times, which, in turn, requires access to all agent observations at all times. Hence, the SPRT may not be implementable beyond the fully centralized or all-to-all agent communication scenario as discussed in Section~\ref{subsec:prel_seq_test}. Motivated by practicable agent networking applications, in this paper we are interested in distributed scenarios, in which inter-agent communication is restricted to a preassigned (possibly sparse) communication graph. In particular, given a graph $G=(V,E)$, possibly sparse, modeling inter-agent communication, we consider scenarios in which inter-agent cooperation is limited to a single round of message exchanges among neighboring agents per observation sampling epoch. To formalize the distributed setup and the corresponding subclass $\mathcal{D}_{d}$ of distributed tests, denote by $\mathcal{G}_{d,i}(t)$ the information available at an agent $i$ at time $t$. The information set includes the observations sampled by $i$ and the messages received from its neighbors till time $t$, and is formally given by the $\sigma$-algebra
\begin{align}
\label{eq:0.02a}
\mathcal{G}_{d,i}(t)=\sigma\left\{y_{i}(s), m_{i,j}(s),~~\forall 1\le s \le t, \forall j\in \Omega_{i}\right\}.
\end{align}
The quantity $m_{i,j}(s)$ denotes the message received by $i$ from its neighbor $j\in\Omega_{i}$ at time $s$, assumed to be a vector of constant (time-invariant) dimension; the exact message generation rule is determined by the (distributed) testing procedure $D_{d}$ in place and, necessarily, $m_{i,j}(s)$ is measurable w.r.t. the $\sigma$-algebra $\mathcal{G}_{d,j}(s)$. Based on the information content $\mathcal{G}_{d,i}(t)$ at time $t$, an agent decides on whether to continue taking observations or to stop in the case of which, it decides on one of the hypothesis $H_{0}$ or $H_{1}$. A distributed testing procedure $D_{d}$ then consists of message generation rules, and, local stopping and decision criteria at the agents. Intuitively, and formally by~\eqref{eq:0.02a} and the fact that $m_{i,j}(s)$ is measurable w.r.t $\mathcal{G}_{d,j}(s)$ for all $(i,j)$ and $s$, we have
\begin{equation}
\label{subsigma}
\mathcal{G}_{d,i}(t)\subset\mathcal{G}_{c}(t)~~~\forall i,t,
\end{equation}
i.e., the information available at an agent $i$ in the distributed setting is a subset of the information that would be available to a hypothetical center in a centralized setting as given in Section~\ref{subsec:prel_seq_test}. Formally, this implies that the class of distributed tests $\mathcal{D}_{d}$ is a subset of the class of centralized or all-possible tests $\mathcal{D}_{c}$ as given in Section~\ref{subsec:prel_seq_test}, i.e., $\mathcal{D}_{d}\subset\mathcal{D}_{c}$. (Intuitively, it means any distributed test can be implemented in a centralized setup or by assuming all-to-all communication.) In this paper, we are interested in characterizing the distributed test that conforms to the communication restrictions above and is optimal in the following sense:
\begin{align}
\label{eq:0.02}
&\min\limits_{D_{d}\in \mathcal{D}_{d}}\max\limits_{i=1,2,\hdots,N}\mathbb{E}_{1}[T_{D_{d},i}], \nonumber\\
&\mbox{s.t.}~\mathbb{P}_{FA}^{D_{d},i}\le\alpha, \mathbb{P}_{M}^{D_{d},i}\le\beta, \forall i=1,2,\hdots,N.
\end{align}
In the above, $T_{D_{d},i}$ denotes the termination (stopping) time at an agent $i$ and $\mathbb{P}_{FA}^{D_{d},i}$, $\mathbb{P}_{M}^{D_{d},i}$, the respective false alarm and miss probabilities at $i$. Note that, since $\mathcal{D}_{d}\subset\mathcal{D}_{c}$, for any distributed test $D_{d}$ we have $\mathbb{E}_{1}[T_{D_{d},i}]\geq \mathbb{E}_{1}[T_{c}]$ for all $i$ at any specified $\alpha$ and $\beta$, i.e., a distributed procedure cannot outperform the optimal centralized procedure, the SPRT given by~\eqref{eq:f5b}-\eqref{eq:f5c}. Rather than solving~\eqref{eq:0.02}, in this paper, we propose a distributed testing procedure of the consensus+innovations type (see Section~\ref{dsd}), which is efficiently implementable and analyze its performance w.r.t. the optimal centralized testing procedure. In particular, we study its performance as a function of the inter-agent communication graph and show that as long as the network is \emph{reasonably} well-connected, but possibly much sparser than the complete or all-to-all network, the suboptimality (in terms of the expected stopping times at the agents) of the proposed distributed procedure w.r.t. the optimal centralized SPRT procedure is upper bounded by a constant factor much smaller than $N$. Our results clearly demonstrate the benefits of collaboration (even over a sparse communication network) as, in contrast, in the non-collaboration case (i.e., each agent relies on its own observations only) each agent would require $N$ times the expected number of observations to achieve prescribed $\alpha$ and $\beta$ as compared to the optimal centralized scenario.\footnote{In the non-collaboration setup, the optimal procedure at an agent is to perform an SPRT using its local observation sequence only; w.r.t. the centralized, this implies that the effective SNR at an agent reduces by a factor $1/N$ and hence (see~\cite{wald1948optimum}) the agent would require $N$ times more observations (in expectation) to achieve the same level of false alarm and miss.}

\section{A Distributed Sequential Detector}
\label{dsd}
\noindent To mitigate the high communication and synchronization overheads in centralized processing, we propose a distributed sequential detection scheme where network communication is restricted to a more \emph{localized} agent-to-agent interaction scenario. More specifically, in contrast to the fully centralized setup described in Section \ref{subsec:prel_seq_test}, we now consider sequential detection in a distributed information setup in which inter-agent information exchange or cooperation is restricted to a preassigned (arbitrary, possibly sparse) communication graph, whereby an agent exchanges its (scalar) test statistic and a scalar function of its latest sensed information with its (one-hop) neighbors.
In order to achieve reasonable detection performance with such localized interaction, we propose a distributed sequential detector of the \emph{consensus}+\emph{innovations} form. Before discussing the details of our algorithm, we state an assumption on the inter-agent communication graph.
\begin{myassump}{A3}
\label{as:3}
\emph{The inter-agent communication graph is connected, i.e. $\lambda_{2}(\mathbf{L}) > 0$, where $\mathbf{L}$ denotes the associated graph Laplacian matrix}.
\end{myassump}

\noindent\textbf{Decision Statistic Update}.
In the proposed distributed algorithm, each agent $i$ maintains a test statistic $P_{d,i}(t)$, which is updated recursively in a distributed fashion as follows :
\begin{align}
\label{eq:2}
&P_{d,i}(t+1)=\frac{t}{t+1}\left(w_{ii}P_{d,i}(t)+\sum_{j\in \Omega_{i}}w_{ij}P_{d,j}(t)\right)\nonumber\\&+\frac{1}{t+1}\left(w_{ii}\eta_{i}(t+1)+\sum_{j\in \Omega_{i}}w_{ij}\eta_{j}(t+1)\right),
\end{align}

\noindent where $\Omega_{i}$ denotes the communication neighborhood of agent $i$ and the $w_{ij}$'s denote appropriately chosen combination weights (to be specified later).

\noindent We collect the weights $w_{ij}$ in an $N\times N$ matrix $\mathbf{W}$, where we assign $w_{ij}=0$, if $(i,j)\notin E$. Denoting by $\mathbf{P}_{d}(t)$ and $\eta(t)$ as the vectors $[P_{d,1}(t),P_{d,2}(t),\hdots,P_{d,N}(t)]^{\top}$ and $[\eta_{1}(t),\eta_{2}(t),\hdots,\eta_{N}(t)]^{\top}$ respectively, \eqref{eq:2} can be compactly written as follows:-
\begin{align}
\label{eq:3}
\mathbf{P_{d}}(t+1)=\mathbf{W}\left(\frac{t}{t+1}\mathbf{P_{d}}(t)+\frac{1}{t+1}\mathbf{\eta}(t+1)\right).
\end{align}

\noindent Now we state some design assumptions on the weight matrix $\mathbf{W}$.
\begin{myassump}{A4}
\label{as:4}
\emph{We design the weights $w_{ij}$'s in \eqref{eq:2} such that the matrix $\mathbf{W}$ is non-negative, symmetric, irreducible and stochastic, i.e., each row of $\mathbf{W}$ sums to one}.
\end{myassump}

\noindent We remark that, if Assumption \ref{as:4} is satisfied, then the second largest eigenvalue in magnitude of $\mathbf{W}$, denoted by $r$, turns out to be strictly less than one, see for example \cite{dimakis2010gossip}.  Note that, by the stochasticity of $\mathbf{W}$, the quantity $r$ satisfies
\begin{align}
\label{eq:3.2}
r=||\mathbf{W}-\mathbf{J}||.
\end{align}

For connected graphs, a simple way to design $\mathbf{W}$ is to assign equal combination weights, in which case we have,
\begin{align}
\label{eq:3.1}
\mathbf{W}=\mathbf{I}-\delta\mathbf{L},
\end{align}
where $\delta$ is a suitably chosen constant. As shown in \cite{xiao2004fast,kar2008sensor}, Assumption \ref{as:4} can be enforced by taking $\delta$ to be in $(0, 2/\lambda_{N}(\mathbf{L}))$. The smallest value of $r$ is obtained by setting $\delta$ to be equal to $2/(\lambda_{2}(\mathbf{L})+\lambda_{N}(\mathbf{L}))$, in which case we have,
\begin{align}
\label{eq:q2}
r=||\mathbf{W}-\mathbf{J}||=\frac{(\lambda_{N}(\mathbf{L})-\lambda_{2}(\mathbf{L}))}{(\lambda_{2}(\mathbf{L})+\lambda_{N}(\mathbf{L}))}.
\end{align}

\begin{Remark}
\label{rm:0}
It is to be noted that Assumption \ref{as:4} can be enforced by appropriately designing the combination weights since the inter-agent communication graph is connected (see Assumption \ref{as:3}). Several weight design techniques satisfying Assumption \ref{as:4} exist in the literature (see, for example, \cite{xiao2004fast}). The quantity $r$ quantifies the rate of information flow in the network, and in general, the smaller the $r$ the faster is the convergence of information dissemination algorithms (such as the consensus or gossip protocol on the graph, see for example \cite{dimakis2010gossip,kar2008sensor,kar2009distributed}). The optimal design of symmetric weight matrices $\mathbf{W}$ for a given network topology that minimizes the value $r$ can be cast as a semi-definite optimization problem \cite{xiao2004fast}.
\end{Remark}

\noindent\textbf{Stopping Criterion for the Decision Update}.
We now provide a stopping criterion for the proposed distributed scheme. To this end, let $S_{d,i}(t)$ denote the quantity $tP_{d,i}(t)$, and let $\gamma_{d,i}^{h}$ and $\gamma_{d,i}^{l}$ be thresholds at an agent $i$ (to be determined later) such that agent $i$ stops and makes a decision only when,
\begin{align}
\label{eq:g1}
 S_{d,i}(t) \notin [\gamma_{d,i}^{l},\gamma_{d,i}^{h}]
\end{align}
for the first time.
The stopping time for reaching a decision at an agent $i$ is then defined as,
\begin{align}
\label{eq:f3}
T_{d,i}=\inf \{t ~~|  S_{d,i}(t) \notin [\gamma_{d,i}^{l},\gamma_{d,i}^{h}]\},
\end{align}
and the following decision rule is adopted at $T_{d,i}$ :
\begin{align}
\label{eq:f5}
H=
\begin{cases}
H_{0} & S_{d,i}(T_{d,i}) \le \gamma_{d,i}^{l}\\
H_{1} & S_{d,i}(T_{d,i}) \ge \gamma_{d,i}^{h}.
\end{cases}
\end{align}
 \\

\noindent We refer to this distributed scheme \eqref{eq:2}, \eqref{eq:f3} and \eqref{eq:f5} as the \emph{consensus+innovations} SPRT ($\cisprt$)  hence forth.

\begin{Remark}
\label{rm:1}
It is to be noted that the decision statistic update rule is distributed and recursive, in that, to realize \eqref{eq:2} each agent needs to communicate its current statistic and a scalar function of its latest sensed observation to its neighbors only; furthermore, the local update rule \eqref{eq:2} is a combination of a consensus term reflecting the weighted combination of neighbors' statistics and a local innovation term reflecting the new sensed information of itself and its neighbors. Note that the stopping times  $T_{d,i}$'s are random and generally take different values for different agents. It is to be noted that the $T_{d,i}$'s are in fact stopping times with respect to the respective agent information filtrations $\mathcal{G}_{d,i}(t)$'s as defined in \eqref{eq:0.02a}. For subsequent analysis we refer to the stopping time of an agent as the stopping time for reaching a decision at an agent.
\end{Remark}

\noindent We end this section by providing some elementary properties of the distributed test statistics.

\begin{proposition}
\label{prop:SdGauss} Let the Assumptions~\ref{as:1}, \ref{as:3} and \ref{as:4} hold. For each $t$ and $i$, the statistic $S_{d,i}(t)$, defined in~\eqref{eq:g1}-\eqref{eq:f5}, is Gaussian under both $H_{0}$ and $H_{1}$. In particular, we have
\begin{align}
\label{prop:SdGauss1}
\mathbb{E}_{0}[S_{d,i}(t)] = -mt~~~\mbox{and}~~~\mathbb{E}_{1}[S_{d,i}(t)]=mt,
\end{align}
and
\begin{align}
\label{prop:SdGauss2}
\mathbb{E}_{0}\left[(S_{d,i}(t)+mt)^{2}\right]=\mathbb{E}_{1}\left[(S_{d,i}(t)-mt)^{2}\right]\leq \frac{2mt}{N}+\frac{2mr^{2}(1-r^{2t})}{1-r^{2}}.
\end{align}
\end{proposition}
\begin{IEEEproof}
Recall from \eqref{eq:f8}, $\eta_{i}(t)$ is distributed as $\mathcal{N}(m,2m)$,~~$\forall i=1,2,\cdots,N$, when conditioned on hypothesis $H_{1}$ and where $m$ is the Kullback-Leibler divergence as defined in \eqref{eq:snr}. Hence,
\begin{align}
\label{eq:13.1}
&\mathbb{E}_{1}[S_{d,i}(t)]=\sum_{j=1}^{t}\mathbf{e}_{i}^{\top}\mathbf{W}^{t+1-j}\mathbb{E}_{1}[\mathbf{\eta}(j)]\nonumber\\
&=m\sum_{j=1}^{t}\mathbf{e}_{i}^{\top}\mathbf{W}^{t+1-j}\mathbf{1}\nonumber\\
&\Rightarrow \mathbb{E}_{1}[S_{d,i}(t)]=mt.
\end{align}
We note that $\mathbf{S_{\eta}}=Cov(\mathbf{\eta}(t))=2m\mathbf{I}$.
By standard algebraic manipulations we have,
\begin{align}
\label{eq:13.2}
&Var(S_{d,i}(t))=\mathbb{E}_{1}\left[(S_{d,i}(t)-mt)^{2}\right]=\sum_{j=1}^{t}\mathbf{e}_{i}^{\top}\mathbf{W}^{t+1-j}\mathbf{S}_{\eta}\mathbf{W}^{t+1-j}\mathbf{e}_{i}\nonumber\\
&=\sum_{j=1}^{t}\mathbf{e_{i}^{\top}(W^{t-j}-J)S_{\eta}(W^{t-j}-J)e_{i}}+\sum_{j=1}^{t}\mathbf{e_{i}^{\top}JS_{\eta}Je_{i}}\nonumber\\
&=2m\sum_{j=1}^{t}\mathbf{e_{i}}^{\top}(\mathbf{W}^{2(t-j)}-\mathbf{J})\mathbf{e_{i}}+2m\sum_{j=1}^{t}\mathbf{e_{i}}^{\top}\mathbf{J}\mathbf{e_{i}}\nonumber\\
&=2m||\sum_{j=1}^{t}\mathbf{e_{i}}^{\top}(\mathbf{W}^{2(t-j)}-\mathbf{J})\mathbf{e_{i}}||+\frac{2mt}{N}\nonumber\\
&\le 2m\sum_{j=1}^{t}||\mathbf{e_{i}}^{\top}(\mathbf{W}^{2(t-j)}-\mathbf{J})\mathbf{e_{i}}||+\frac{2mt}{N}\nonumber\\
&\le 2m\sum_{j=1}^{t}||\mathbf{e_{i}}^{T}\mathbf{e_{i}}|| ||\mathbf{W}^{2(t-j)}-\mathbf{J}||+\frac{2mt}{N}\nonumber\\
&=2m\sum_{j=0}^{t-1}r^{2j}++\frac{2mt}{N}\nonumber\\
&\le\frac{2mt}{N}+\frac{2m(1-r^{2t})}{1-r^{2}}.
\end{align}
The assertion for hypothesis $H_{0}$ follows in a similar way.
\end{IEEEproof}

\section{Main Results}
\label{sec:main_res}
\noindent We formally state the main results in this section, the proofs being provided in Section~\ref{sec:proof_res}.

\subsection{Thresholds for the $\cisprt$}
\noindent In this section we derive thresholds for the $\cisprt$, see \eqref{eq:g1}-\eqref{eq:f5}, in order to ensure that the procedure terminates in finite time a.s. at each agent and the agents achieve specified error probability requirements.
We emphasize that in the proposed approach, a particular agent has access to its one hop neighborhood's test statistics and latest sensed information only. Moreover the latest sensed information is accessed through a scalar function of the latest observation of the agents in an agent's neighborhood. Recall, by \eqref{eq:2} and \eqref{eq:g1} the (distributed) test statistic at agent $i$ is given by
\begin{align}
\label{eq:dist_stat}
S_{d,i}(t)=\sum_{j=1}^{t}\mathbf{e}_{i}^{\top}\mathbf{W}^{t+1-j}\mathbf{\eta}_{j}.
\end{align}

\noindent For the proposed  $\cisprt$, we intend to derive thresholds  which guarantee the error performance in terms of the error probability requirements $\alpha$ and $\beta$, i.e., such that $\mathbb{P}_{FA}^{d,i}\le \alpha$ and $\mathbb{P}_{M}^{d,i}\le \beta, ~~\forall i=1, 2, \hdots, N$, where $\mathbb{P}_{FA}^{d,i}$ and $\mathbb{P}_{M}^{d,i}$ represent the probability of false alarm and the probability of miss for the $i$th agent defined as
\begin{align}
\label{eq:dth1}
&\mathbb{P}_{FA}^{d,i}=\mathbb{P}_{0}(S_{d,i}(T_{d,i})\ge\gamma_{d,i}^{h})\nonumber\\
&\mathbb{P}_{M}^{d,i}=\mathbb{P}_{1}(S_{d,i}(T_{d,i})\le\gamma_{d,i}^{l}),
\end{align}
with $T_{d,i}$ as defined in \eqref{eq:f3}.

\begin{Theorem}
\label{dist_th}
Let the Assumptions~\ref{as:1}-\ref{as:4} hold.\\
1) Then, for each $\alpha$ and $\beta$ there exist $\gamma_{d,i}^{h}$ and $\gamma_{d,i}^{l}$,~$\forall i=1,2,\hdots,N$, such that $\mathbb{P}_{FA}^{d,i} \le \alpha$ and $\mathbb{P}_{M}^{d,i} \le \beta$ and the test concludes in finite time a.s. i.e.\\
\begin{align}
\label{eq:th7.1}
\mathbb{P}_{1}(T_{d,i}<\infty)= 1, \forall i=1,2,\hdots,N ,
\end{align}
where $T_{d,i}$ is the stopping time for reaching a decision at agent $i$.\\
2) In particular, for given $\alpha$ and $\beta$, any choice of thresholds $\gamma_{d,i}^{h}$ and $\gamma_{d,i}^{l}$ satisfying
\begin{align}
\label{eq:d1}
\gamma_{d,i}^{h}\ge\frac{8(k+1)}{7N}\left(\log\left(\frac{2}{\alpha}\right)-\log(1-e^{\frac{-Nm}{4(k+1)}})\right)=\gamma_{d}^{h,0}
\end{align}
\begin{align}
\label{eq:d2}
\gamma_{d,i}^{l}\le\frac{8(k+1)}{7N}\left(\log\left(\frac{\beta}{2}\right)+\log(1-e^{\frac{-Nm}{4(k+1)}})\right)=\gamma_{d}^{l,0},
\end{align}
where $m$ is defined in \eqref{eq:snr} and $k$ is defined by
\begin{align}
\label{eq:d3}
Nr^{2}= k,
\end{align}
with $r$ as in \eqref{eq:3.2}, achieves a.s. finite stopping at an agent $i$ while ensuring that $\mathbb{P}_{FA}^{d,i} \le \alpha$ and $\mathbb{P}_{M}^{d,i} \le \beta$.
\end{Theorem}

\noindent The first assertion ensures that for any set of pre-specified error metrics $\alpha$ and $\beta$ (satisfying Assumption \ref{as:2}), the $\cisprt$ can be designed to achieve the error requirements while ensuring finite stopping a.s.
It is to be noted that the ranges associated with the thresholds in \eqref{eq:d1}-\eqref{eq:d2} provide sufficient threshold design conditions for achieving pre-specified performance, but may not be necessary. The thresholds chosen according to \eqref{eq:d1}-\eqref{eq:d2}  are not guaranteed to be optimal in the sense of the expected stopping time of the $\cisprt$ algorithm and there might exist better thresholds (in the sense of expected stopping time) that achieve the pre-specified error requirements.

\begin{Remark}
\label{rm:31}
We remark the following: 1) We have shown that the $\cisprt$ algorithm can be designed so as to achieve the pre-specified error metrics at every agent $i$. This, in turn, implies that the probability of not reaching decision consensus among the agents can be upper bounded by $N\beta$ when conditioned on $H_{1}$ and $N\alpha$ when conditioned on $H_{0}$. It is to be noted that with $\alpha \to 0$ and $\beta \to 0$, the probability of not reaching decision consensus conditioned on either of the hypothesis goes to $0$ as well;
2) The factor $k$ in the closed form expressions of the thresholds in \eqref{eq:d1} and \eqref{eq:d2} relates the value of the thresholds to the rate of flow of  information $r$ and, hence, in turn, can be related to the degree of connectivity of the inter-agent communication graph under consideration, see \eqref{eq:3.2}-\eqref{eq:3.1} and the accompanying discussion. From Assumption \ref{as:4}, we have that $r < 1$. As $r$ goes smaller, which intuitively means increased rate of flow of information in the inter-agent network, the value of thresholds needed to achieve the pre-specified error metrics become smaller i.e. the interval $[\gamma_{d,i}^{l}, \gamma_{d,i}^{h}]$ shrinks for all $i=1, 2, \hdots, N$.
\end{Remark}

\subsection{Probability Distribution of $T_{d,i}$ and $T_{c}$}

\noindent We first characterize the stopping time distributions for the centralized SPRT detector (see Section \ref{subsec:prel_seq_test}) and those of the distributed $\cisprt$. Subsequently, we compare the centralized and distributed stopping times by studying their respective large deviation tail probability decay rates.
\begin{Theorem}(\hspace{-0.5pt}\cite{darling1953first,hieber2012note})
\label{1.1}
Let the Assumptions \ref{as:1} and \ref{as:2} hold and
given the SPRT for the centralized setup in \eqref{eq:cent_stat}-\eqref{eq:f5c}, we have
\begin{align}
\label{eq:1.01}
\mathbb{P}_{1}(T_{c}> t)\geq\exp\left(\frac{N\mu\gamma_{c}^{l}}{\sigma^{2}}\right)K_{t}^{\infty}\left(\gamma_{c}^{h}\right)-\exp\left(\frac{N\mu\gamma_{c}^{h}}{\sigma^{2}}\right)K_{t}^{\infty}\left(\gamma_{c}^{l}\right),
\end{align}
where
\begin{align}
\label{eq:1.011}
K_{t}^{S}\left(a\right)=\frac{\sigma^{2}\pi}{N(\gamma_{c}^{h}-\gamma_{c}^{l})^{2}}\sum_{s=1}^{S}\frac{l(-1)^{l+1}}{\frac{Nm}{4}+\frac{\sigma^{2}s^{2}\pi^{2}}{2N(\gamma_{c}^{h}-\gamma_{c}^{l})^{2}}}\exp\left(-\left(\frac{Nm}{4}+\frac{\sigma^{2}s^{2}\pi^{2}}{2N(\gamma_{c}^{h}-\gamma_{c}^{l})^{2}}\right)t\right)\sin\left(\frac{s\pi a}{\gamma_{c}^{h}-\gamma_{c}^{l}}\right),
\end{align}

whereas, $T_{c}$ is defined in \eqref{eq:f5b} and $\gamma_{c}^{h}$ and $\gamma_{c}^{l}$ are the associated SPRT thresholds chosen to achieve specified error requirements $\alpha$ and $\beta$.
\end{Theorem}

\noindent The above characterization of the stopping distribution of Wald's SPRT was obtained in \cite{darling1953first,hieber2012note}. In particular, this was derived by studying the first passage time distribution of an associated continuous time Wiener process with a constant drift; intuitively, the continuous time approximation of the discrete time SPRT consists of replacing the discrete time likelihood increments by a Wiener process accompanied by a constant drift that reflects the mean of the hypothesis in place. This way, the sequence obtained by sampling the continuous time process at integer time instants is equivalent in distribution to the (discrete time) Wald's SPRT considered in this paper. The term on the R.H.S. of~\eqref{eq:1.01} is exactly equal to the probability that the first passage time of the continuous time Wiener process with left and right boundaries $\gamma_{c}^{l}$ and $\gamma_{c}^{h}$ respectively is greater than $t$, whereas, is, in general, a lower bound for the discrete time SPRT (as given in Theorem~\ref{1.1}) as increments in the latter happen at discrete (integer) time instants only.

\noindent We now provide a characterization of the stopping time distributions of the $\cisprt$ algorithm.

\begin{Lemma}
\label{1.2}
Let the assumptions \ref{as:1}-\ref{as:4} hold. Consider the  $\cisprt$ algorithm given in \eqref{eq:3}, \eqref{eq:f3} and \eqref{eq:f5} and suppose that, for specified $\alpha$ and $\beta$, the thresholds $\gamma_{d,i}^{h}$ and $\gamma_{d,i}^{l}$, $i=1,\cdots,N$, are chosen to satisfy the condtions derived in \eqref{eq:d1} and \eqref{eq:d2}. We then have,
\begin{align}
\label{eq:1.02}
&\mathbb{P}_{1}(T_{d,i}>t) \le\mathbb{Q}\Big(\frac{-\gamma_{d,i}^{h}+mt}{\sqrt{\frac{2mt}{N}+\frac{2mr^{2}(1-r^{2t})}{1-r^{2}}}}\Big),~~\forall i=1, 2, \hdots, N,
\end{align}
where $T_{d,i}$ is the stopping time of the $i$-th agent to reach a decision as defined in~\eqref{eq:f3}.
\end{Lemma}

\subsection{Comparison of stopping times of the distributed and centralized detectors}
\noindent In this section we compare the stopping times $T_{c}$ and $T_{d,i}$ by studying their respective large deviation tail probability decay rates. We utilize the bounds derived in Theorem \ref{1.1} and Lemma \ref{1.2} to this end.
\begin{corollary}
\label{7.2}
Let the hypotheses of Lemma \ref{1.1} hold. Then we have the following large deviation characterization for the tail probabilities of $T_{c}$:
\begin{align}
\label{eq:12.52}
\liminf_{t\to\infty}\frac{1}{t}\log(\mathbb{P}_{1}(T_{c}> t))\geq -\frac{Nm}{4}-\frac{\sigma^{2}\pi^{2}}{2N(\gamma_{c}^{h}-\gamma^{l})^{2}}.
\end{align}
\end{corollary}

It is to be noted that the exponent is a function of the thresholds $\gamma_{c}^{h}$ and $\gamma_{c}^{l}$ and with the decrease in the error constraints $\alpha$ and $\beta$,  $\frac{Nm}{4}+\frac{\sigma^{2}\pi^{2}}{2N(\gamma_{c}^{h}-\gamma^{l})^{2}} \approx \frac{Nm}{4}$.

\begin{Theorem}
\label{7.3}
Let the hypotheses of Lemma \ref{1.2} hold. Then we have the following large deviation characterization for the tail probabilities of the $T_{d,i}$'s:
\begin{align}
\label{eq:12.69}
\limsup_{t\to\infty}\frac{1}{t}\log(\mathbb{P}_{1}(T_{d,i}> t))\le -\frac{Nm}{4} ~~,\forall i=1,2,\hdots,N.
\end{align}
\end{Theorem}

\noindent
Importantly, the upper bound for the large deviation exponent of the $\cisprt$ in Theorem \ref{7.3} is independent of the inter-agent communication topology as long as the connectivity conditions Assumptions \ref{as:3}-\ref{as:4} hold. Finally, in the asymptotic regime, i.e., as $N$  goes to $\infty$, since $\frac{\sigma^{2}\pi^{2}}{2N(\gamma_{c}^{h}-\gamma^{l})^{2}}=o(Nm)$, we have that the performance of the distributed $\cisprt$ approaches that of the centralized SPRT, in the sense of stopping time tail exponents, as $N$ tends to $\infty$.

\subsection{Comparison of the expected stopping times of the centralized and distributed detectors}
\noindent In this section we compare the expected stopping times of the centralized SPRT detector and the  proposed $\cisprt$ detector.  Recall that $\mathbb{E}_{j}[T_{d,i}]$ and $\mathbb{E}_{j}[T_{c}]$ represent the expected stopping times for reaching a decision for the $\cisprt$ (at an agent $i$) and its centralized counterpart respectively, where $j\in\{0,1\}$ denotes the hypothesis on which the expectations are conditioned on. Without loss of generality we compare the expectations conditioned on Hypothesis $H_{1}$, similar conclusions (with obvious modifications) hold when the expectations are conditioned on $H_{0}$ (see also Section \ref{subsec:prel_seq_test}).

\noindent Also, for the sake of mathematical brevity and clarity, we approximate $\alpha=\beta=\epsilon$ in this subsection.

\noindent Recall Section \ref{subsec:prel_seq_test} and note that, at any instant of time $t$, the information $\sigma$-algebra $\mathcal{G}_{d,i}(t)$ at any agent $i$ is a subset of $\mathcal{G}_{c}(t)$, the information $\sigma$-algebra of a (hypothetical) center, which has access to the data of all agents at all times. This implies that any distributed procedure (in particular the $\cisprt$) can be implemented in the centralized setting, and, since $M(\epsilon)$ (see~\eqref{eq:edc3100}) constitutes a lower bound on the expected stopping time of any sequential test achieving error probabilities $\alpha=\beta=\epsilon$, we have that
\begin{align}
\label{eq:edc1}
\frac{\mathbb{E}_{1}[T_{d,i}]}{\mathcal{M}(\epsilon)}\ge 1,~~\forall~~i=1,2,\hdots,N,
\end{align}
for all $\epsilon\in (0,1/2)$.
\noindent In order to provide an upper bound on the ratio $\mathbb{E}_{1}[T_{d,i}]/\mathcal{M}(\epsilon)$ and, hence, compare the performance of the proposed $\cisprt$ detector with the optimal centralized detector, we first obtain a characterization of $\mathbb{E}_{1}[T_{d,i}]$ in terms of the algorithm thresholds as follows.
\begin{Theorem}
\label{edct_1}
Let the assumptions \ref{as:1}-\ref{as:4} hold and let $\alpha=\beta=\epsilon$. Suppose that the thresholds of the $\cisprt$ be chosen as $\gamma_{d,i}^{h}=\gamma_{d}^{h,0}$ and $\gamma_{d,i}^{l}=\gamma_{d}^{l,0}$ for all $i=1,\cdots,N$, where $\gamma_{d}^{h,0}$ and $\gamma_{d}^{l,0}$ are defined in \eqref{eq:d1}-\eqref{eq:d2}. Then, the stopping time $T_{d,i}$ of the $\cisprt$ at an agent $i$ satisfies
\begin{align}
\label{eq:edct_11}
\frac{(1-2\epsilon)\gamma_{d,i}^{h}}{m}-\frac{c}{m}\le\mathbb{E}_{1}[T_{d,i}] \le \frac{5\gamma_{d,i}^{h}}{4m}+\frac{1}{1-e^{\frac{-Nm}{4(k+1)}}},
\end{align}
where $k=Nr^{2}$, $r$ is as defined in \eqref{eq:3.2}, and $c>0$ is a constant that may be chosen to be independent of the thresholds and the $\epsilon$.
\end{Theorem}

\noindent It is to be noted that, when $\alpha=\beta=\epsilon$, then $\gamma_{d,i}^{h}=-\gamma_{d,i}^{l}$ from \eqref{eq:d1} and \eqref{eq:d2}. The upper bound derived in the above assertion might be loose, owing to the approximations related to the non-elementary $\mathbb{Q}$-function. We use the derived upper bound for comparing the performance of the $\cisprt$ algorithm with that of its centralized counterpart. The constant $c>0$ in the lower bound is independent of the thresholds $\gamma_{d,i}^{l}$ and $\gamma_{d,i}^{h}$ (and hence, also independent of the error tolerance $\epsilon$) and is a function of the network topology and the Gaussian model statistics only. Explicit expressions and bounds on $c$ may be obtained by refining the various estimates in the proofs of Lemma~\ref{lm:div_est} and Theorem~\ref{edct_1}, see Section~\ref{sec:main_res}. However, for the current purposes, it is important to note that $c=o(\gamma_{d}^{h,0})$, i.e., as $\epsilon$ goes to zero or equivalently in the limit of large thresholds $c/\gamma_{d}^{h,0}\rightarrow 0$. Hence, as $\epsilon\rightarrow 0$, the more readily computable quantity $\frac{(1-2\epsilon)\gamma_{d,i}^{h}}{m}$ may be viewed as a reasonably good approximation to the lower bound in Theorem~\ref{edct_1}.

\begin{Theorem}
\label{edct}
Let the hypotheses of Theorem \ref{edct_1} hold. Then, we have the following characterization of the ratio of the expected stopping times of the $\cisprt$ and the centralized detector in asymptotics of the $\epsilon$,
\begin{align}
\label{eq:edc4}
1\le\limsup_{\epsilon\to0}\frac{\mathbb{E}_{1}[T_{d,i}]}{\mathcal{M}(\epsilon)}\le\frac{10(k+1)}{7},~~\forall i=1,2,\hdots,N,
\end{align}
where $k=Nr^{2}$ and $r$ is as defined in \eqref{eq:3.2}.
\end{Theorem}

\noindent Theorem \ref{edct} shows that the $\cisprt$ algorithm can be designed in such a way that with pre-specified error metrics $\alpha$ and $\beta$ going to $0$ , the ratio of the expected stopping time for the $\cisprt$ algorithm and its centralized counterpart are bounded above by $\frac{10(k+1)}{7}$ where the quantity $k$ depends on $r$ which essentially quantifies the dependence of the $\cisprt$ algorithm on the network connectivity.

\begin{Remark}
\label{rm:2}
It is to be noted that the derived upper bound for the ratio of the expected stopping times of the $\cisprt$ algorithm and its centralized counterpart may not be a tight upper bound. The looseness in the upper bound is due to the fact that the set of thresholds chosen are oriented to be sufficient conditions and not necessary. As pointed out in Remark \ref{rm:31} there might exist possibly better choice of thresholds for which the pre-specified error metrics are satisfied. Hence, given a set of pre-specified error metrics and a network topology the upper bound of the derived assertion above can be minimized by choosing the optimal weights for $\mathbf{W}$ as shown in \cite{xiao2004fast}. It can be seen that the ratio of expected stopping times of the isolated SPRT based detector case, i.e., the non-collaboration case, and the centralized SPRT based detector is $N$ (see Section~\ref{subsec:prel_seq_test}). So, for the $\cisprt$ case in order to make savings as far as the stopping time is concerned with respect to the isolated SPRT based detector, $\frac{10(k+1)}{7}\le N$  should be satisfied. Hence, we have that $r\le \sqrt{\frac{7N-10}{10N}}$ is a sufficient condition for the same.
\end{Remark}

\section{Dependence of the $\cisprt$ on Network Connectivity: Illustration}
\label{illust}

\noindent In this section, we illustrate the dependence of the $\cisprt$ algorithm on the network connectivity, by considering a class of graphs. Recall from section \ref{dsd} that  the quantity $r$ quantifies the \emph{rate of information flow} in the network, and in general, the smaller the $r$ the faster is the convergence of information dissemination algorithms (such as the consensus or gossip protocol (\hspace{-0.5pt}\cite{dimakis2010gossip,kar2008sensor,kar2009distributed}) on the graph and the optimal design of symmetric weight matrices $\mathbf{W}$ for a given network topology that minimizes the value $r$ can be cast as a semi-definite optimization problem \cite{xiao2004fast}.

To quantify the dependance of the $\cisprt$ algorithm on the graph topology, we note that the limit derived in \eqref{eq:edc4} is a function of $\mathbf{W}$ and can be re-written as follows :
\begin{align}
\label{eq:i1}
\limsup_{\epsilon\to0}\frac{\mathbb{E}_{1}[T_{d,i}]}{\mathcal{M}(\epsilon)}\le\frac{10(Nr^{2}+1)}{7}\doteq\mathcal{R}(\mathbf{W}),
\end{align}
i.e., the derived upper bound $\mathcal{R}(\mathbf{W})$ is a function of the chosen weight matrix $W$. Based on~\eqref{eq:i1}, naturally, a weight design guideline would be to design $\mathbf{W}$ (under the network topological constraints) so as to minimize $\mathcal{R}(\mathbf{W})$, which, by~\eqref{eq:i1} and as discussed earlier corresponds to minimizing $r=\|\mathbf{W}-\mathbf{J}\|$. This leads to the following upper bound on the achievable performance of the $\cisprt$:
\begin{equation}
\label{23456}
\limsup_{\epsilon\to0}\frac{\mathbb{E}_{1}[T_{d,i}]}{\mathcal{M}(\epsilon)}\le\min_{\mathbf{W}}\mathcal{R}(\mathbf{W}).
\end{equation}
By restricting attention to constant link weights, i.e., $\mathbf{W}$'s of the form $(\mathbf{I}-\delta\mathbf{L})$ and noting that
\begin{equation}
\label{234567}
\min_{\delta}\|\mathbf{I}-\delta\mathbf{L}-\mathbf{J}\|=\frac{(\lambda_{N}(\mathbf{L})-\lambda_{2}(\mathbf{L}))}{(\lambda_{2}(\mathbf{L})+\lambda_{N}(\mathbf{L}))},
\end{equation}
(see~\eqref{eq:q2}), we further obtain
\begin{align}
\label{eq:i3}
\limsup_{\epsilon\to0}\frac{\mathbb{E}_{1}[T_{d,i}]}{\mathcal{M}(\epsilon)}\le \min_{\mathbf{W}}\mathcal{R}(\mathbf{W})\le \min_{\delta}\mathcal{R}(\mathbf{I}-\delta\mathbf{L}) =\frac{10}{7}+\frac{10N(\lambda_{N}(\mathbf{L})-\lambda_{2}(\mathbf{L}))^{2}}{7(\lambda_{2}(\mathbf{L})+\lambda_{N}(\mathbf{L}))^{2}}.
\end{align}
The final bound obtained in~\eqref{eq:i3} might not be tight, being an upper bound (there may exist $\mathbf{W}$ matrices not of the form $\mathbf{I}-\delta\mathbf{L}$ with smaller $r$) to a possibly loose upper bound derived in \eqref{eq:edc4}, but, nonetheless, directly relates the performance of the $\cisprt$ to the spectra of the graph Laplacian and hence the graph topology. From~\eqref{eq:i3} we may further conclude that networks with smaller value of the ratio $\lambda_{2}(\mathbf{L})/\lambda_{N}(\mathbf{L})$ tend to achieve better performance. This leads to an interesting graph design question: given resource constraints, specifically, say a restriction on the number of edges of the graph, how to design inter-agent communication networks that tend to minimize the eigen-ratio $\lambda_{2}(\mathbf{L})/\lambda_{N}(\mathbf{L})$ so as to achieve improved $\cisprt$ performance. To an extent, such graph design questions have been studied in prior work, see~\cite{kar2008topology}, which, for instance, shows that expander graphs tend to achieve smaller $\lambda_{2}(\mathbf{L})/\lambda_{N}(\mathbf{L})$ ratios given a constraint on the total number of network edges.

\section{Simulations}
\label{sec:sim}
\begin{figure}
	\centering
	\captionsetup{justification=centering}
	\includegraphics[width=90mm]{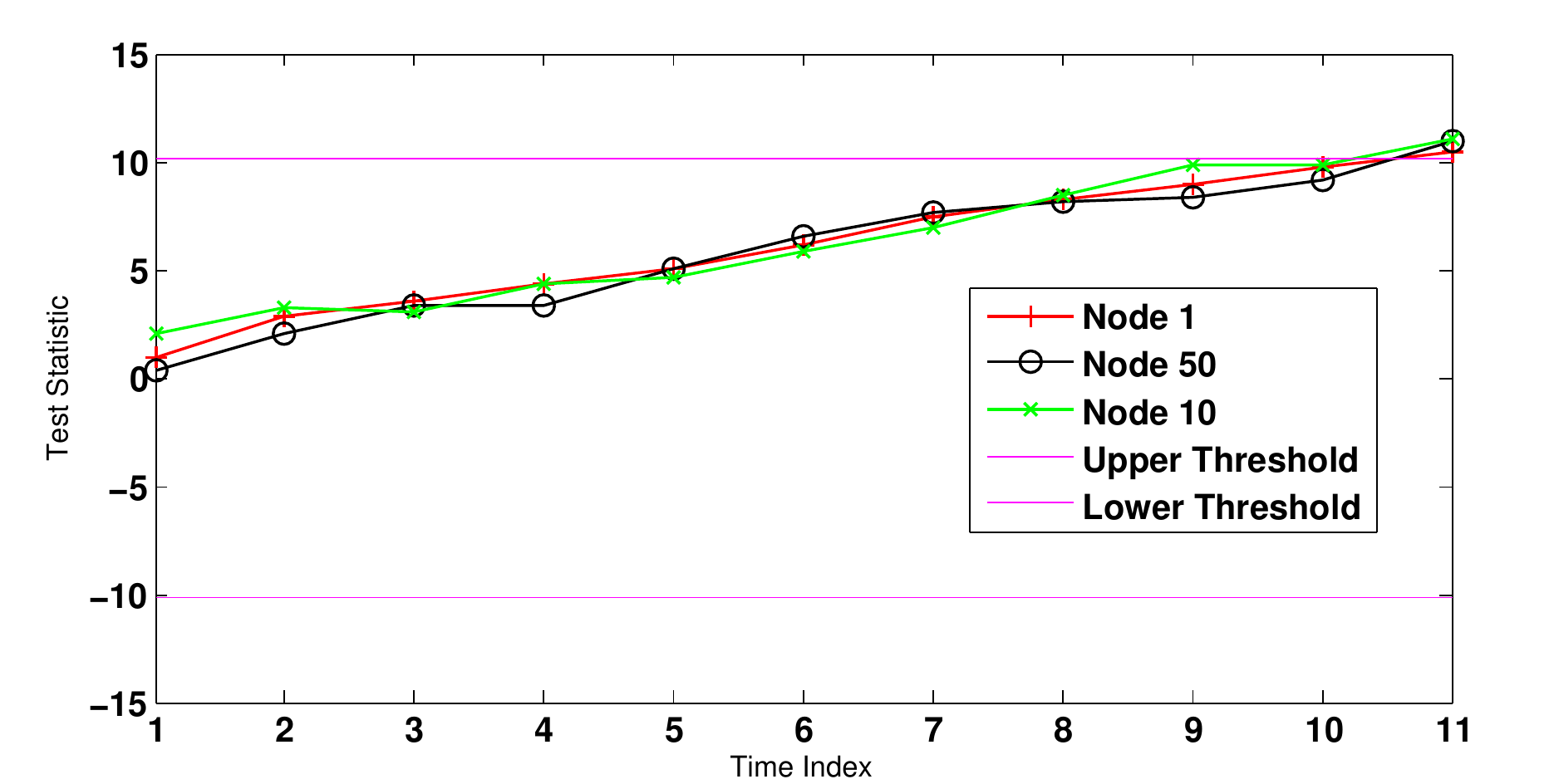}
	\caption{Instantaneous behavior of $S_{d,i}(t)$}\label{fig:4}
\end{figure}

We generate planar random geometric networks of $30$, $300$ and $1000$ agents. The $x$ coordinates and the $y$ coordinates of the agents are sampled from an uniform distribution on the open interval $(0,1)$. We link two vertices by an edge if the distance between them is less than or equal to $g$. We go on re-iterating this procedure until we get a connected graph. We construct the geometric network for each of $N=30, 300~\textrm{and}~1000$ cases with three different values of $g$ i.e. $g=0.3, 0.6~\textrm{and}~0.9$. The values of $r$ obtained in each case is specified in Table \ref{tab:1}.
\begin{table}[h]
	\centering
	\begin{tabular}{| l | l | l | l |}
		\hline
		r & g=0.3 & g=0.6 & g=0.9 \\ \hline
		N=30 & 0.8241 & 0.5580 & 0.2891 \\ \hline
		N=300 & 0.7989 & 0.6014 & 0.2166 \\ \hline
		N=1000 & 0.7689 & 0.5940 & 0.2297 \\
		\hline
	\end{tabular}
	\caption{Values of $r$}\label{tab:1}
\end{table}
We consider two cases, the $\cisprt$ case and the non-collaborative case. We consider $\alpha=\beta=\epsilon$ and ranging from $10^{-8}$ to $10^{-4}$ in steps of $10^{-6}$. For each such $\epsilon$, we conduct $2000$ simulation runs to empirically estimate the stopping time distribution $\mathbb{P}_{1}( T > t)$ of a randomly chosen agent (with uniform selection probability) for each of the cases. From these empirical probability distributions of the stopping times, we estimate the corresponding expected stopping times. Figure \ref{fig:4} shows the instantaneous behavior of the test statistics in the case of $N=300$ with $\epsilon=10^{-10}$.
\begin{figure}
	\centering
	\captionsetup{justification=centering}
	\includegraphics[width=90mm]{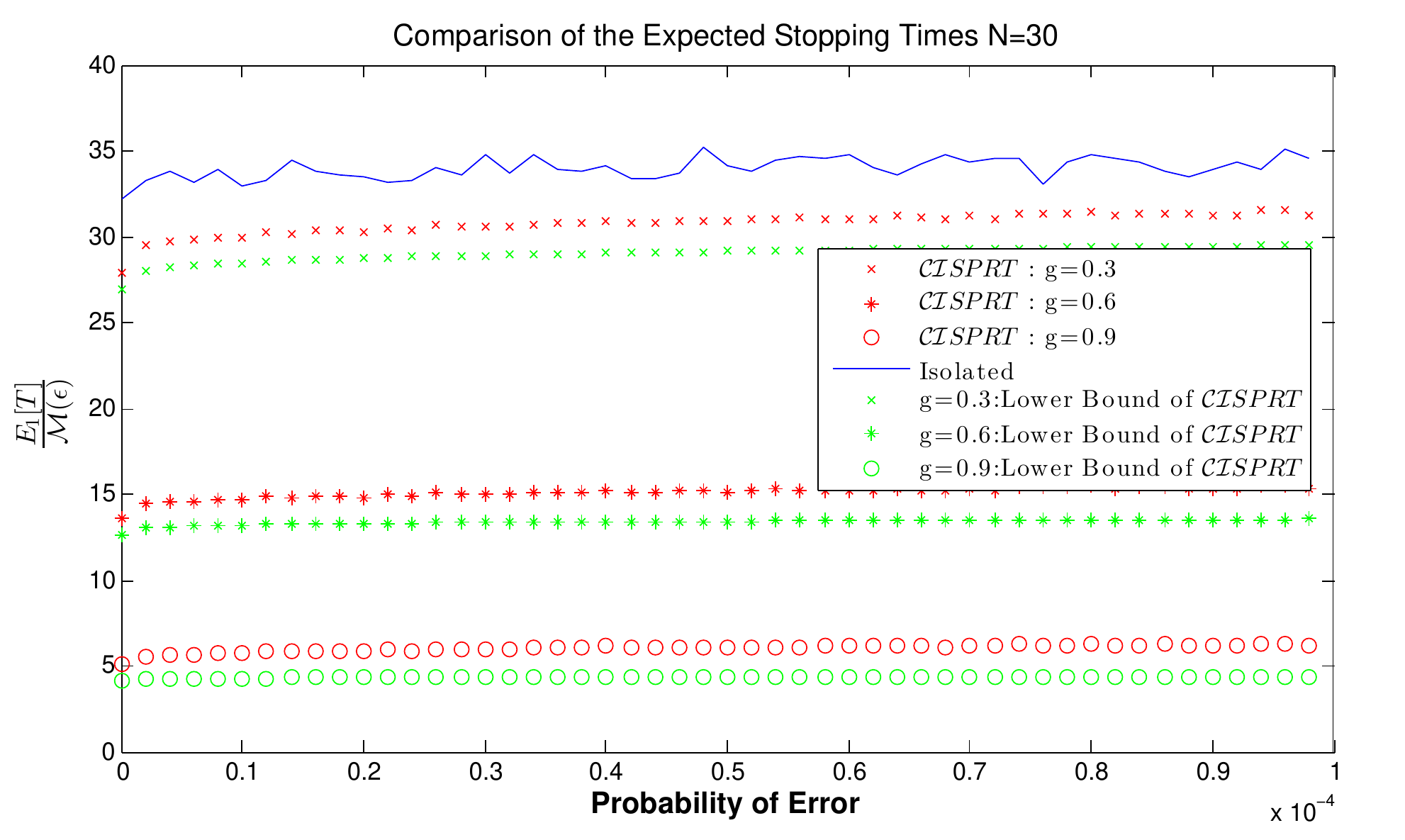}
	\caption{Comparison of Stopping Time Distributions for N=30}\label{fig:1}
\end{figure}
\begin{figure}
	\centering
	\captionsetup{justification=centering}
	\includegraphics[width=90mm]{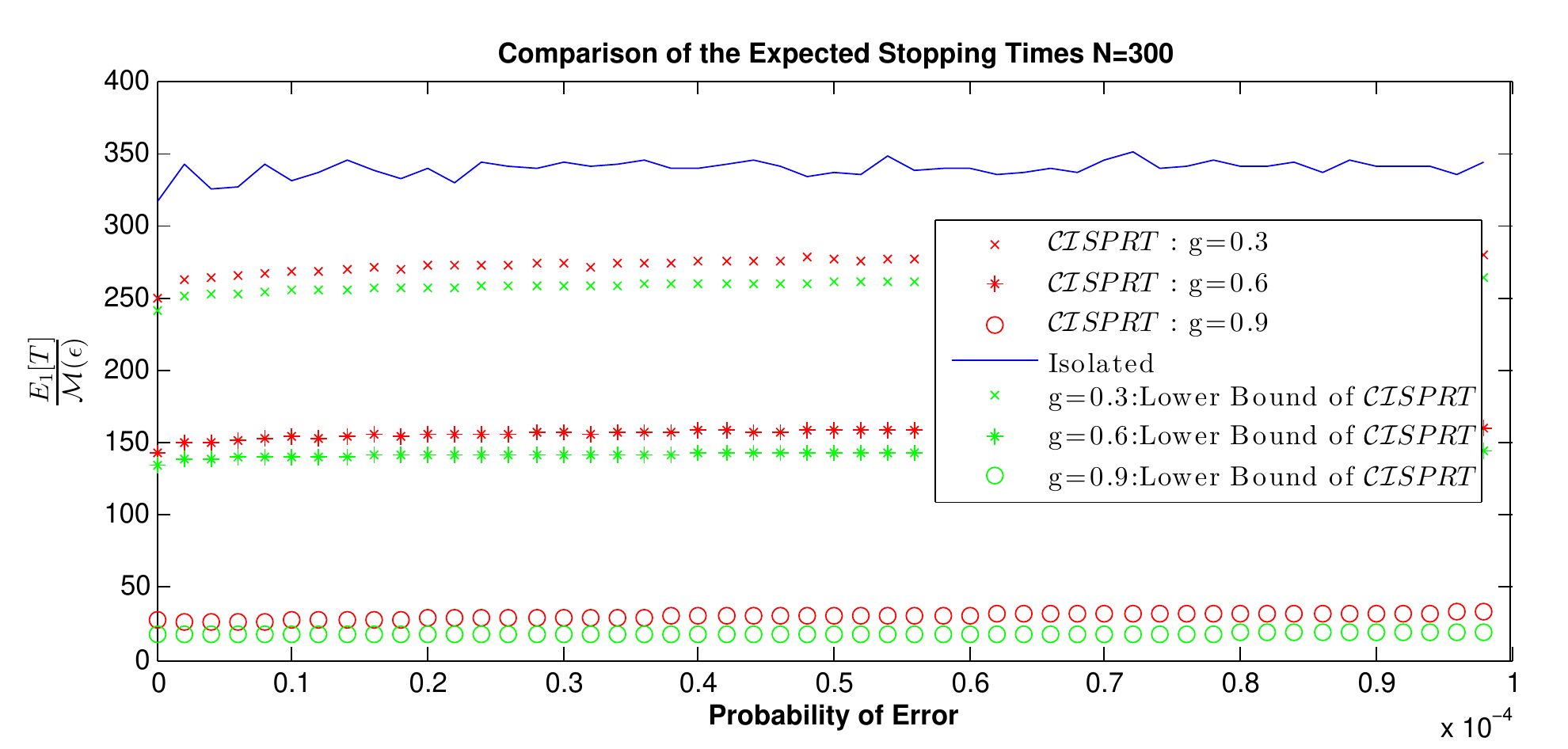}
	\caption{Comparison of Stopping Time Distributions for N=300}\label{fig:2}
\end{figure}
\begin{figure}
	\centering
	\captionsetup{justification=centering}
	\includegraphics[width=90mm]{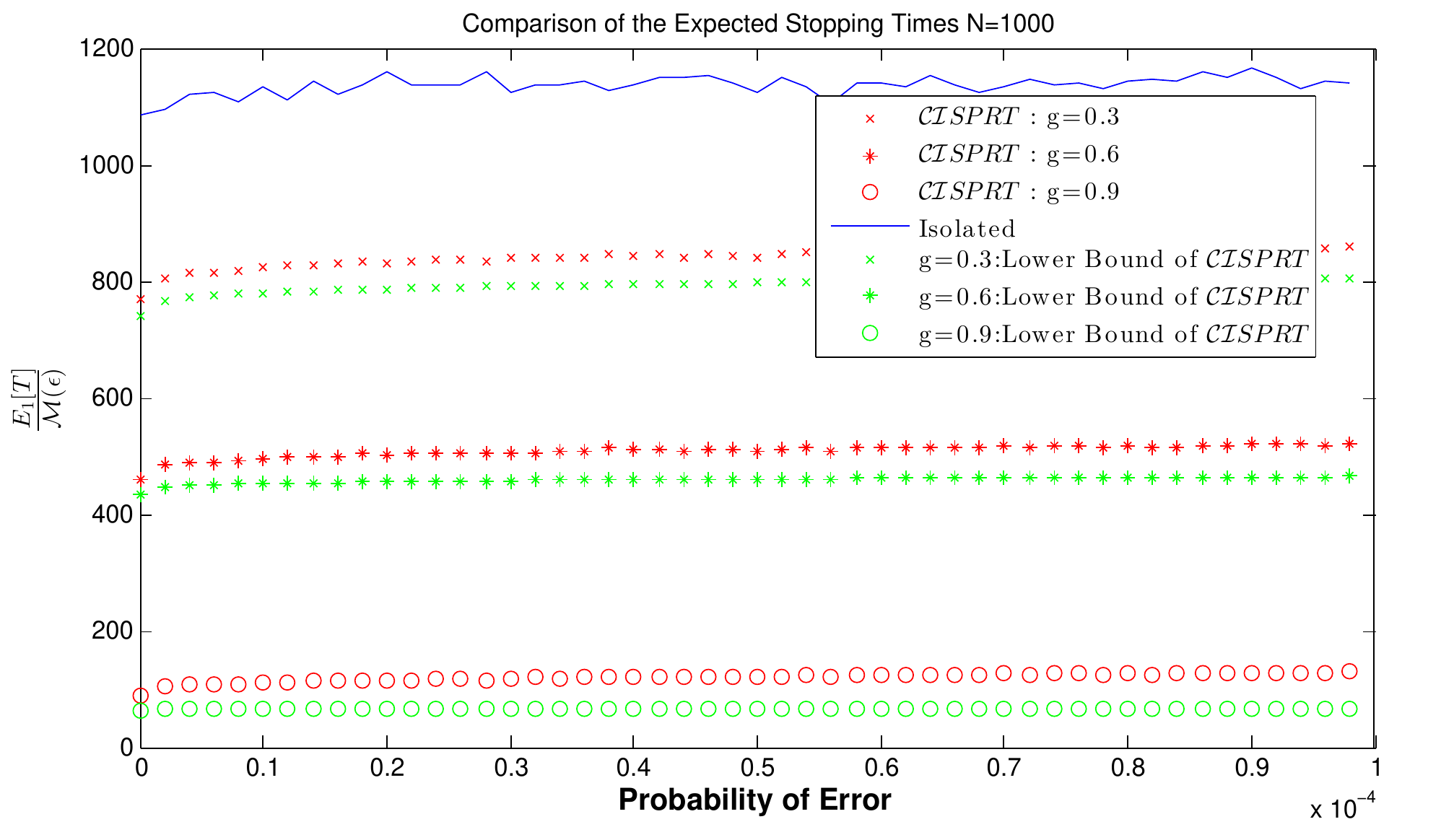}
	\caption{Comparison of Stopping Time Distributions for N=1000}\label{fig:3}
\end{figure}
In Figures \ref{fig:1}, \ref{fig:2} and \ref{fig:3} it is demonstrated that the ratio of the expected stopping time of the $\cisprt$ algorithm and the universal lower bound $\mathcal{M}(\epsilon)$ is less than that of the ratio of the expected stopping times of the isolated (non-collaborative) case and $\mathcal{M}(\epsilon)$. The ratio of the theoretical lower bound of the expected stopping time of the $\cisprt$ derived in Theorem \ref{edct_1} and  $\mathcal{M}(\epsilon)$ was also studied. More precisely, we compared the experimental ratio of the expected stopping times of the $\cisprt$ and $\mathcal{M}(\epsilon)$ with the ratio of the quantity $\frac{(1-2\epsilon)\gamma_{d,i}^{h}}{m}$ (the small $\epsilon$ approximation of the theoretical lower bound given in Theorem~\ref{edct_1}, see also the discussion provided in Section~\ref{sec:main_res} after the statement of Theorem~\ref{edct_1}) and $\mathcal{M}(\epsilon)$. It can be seen that the experimental ratio of the expected stopping times of the $\cisprt$ and $\mathcal{M}(\epsilon)$ is very close to the ratio of the (approximate) theoretical lower bound of expected stopping time of the $\cisprt$ and  $\mathcal{M}(\epsilon)$, which shows that the lower bound derived in Theorem \ref{edct_1} is reasonable. Figure \ref{fig:4} is an example of a single run of the algorithm which shows the instantaneous behavior of the distributed test statistic $S_{d,i}(t)$ for $N=300$, where we have plotted three randomly chosen agents i.e. $i=1$, $i=10$ and $i=50$.

\section{Proofs of Main Results}
\label{sec:proof_res}

\noindent
\begin{IEEEproof}[Proof of Theorem \ref{dist_th}]
Let $\hat{A}=e^{\gamma_{d,i}^{l}}$ and $\hat{B}=e^{\gamma_{d,i}^{h}}$ where $\gamma_{d,i}^{h}$ and $\gamma_{d,i}^{l}$ $ \in \mathbb{R}$ are thresholds (to be designed) for the $\cisprt$.
In the following derivation, for a given random variable $z$ and an event $A$, we use the notation $\mathbb{E}[z;A]$ to denote the expectation $\mathbb{E}[z\mathbb{I}_{A}]$.
Let $T$ denote the random time which can take values in $\overline{\mathbb{Z}}_{+}$ given by
\begin{align}
\label{def:T}
T = \inf\left\{t | S_{d,i}(t) \notin \left[\gamma_{d,i}^{l}, \gamma_{d,i}^{h}\right]\right\}.
\end{align}
First, we show that for any $\gamma_{d,i}^{h}$ and $\gamma_{d,i}^{l}$ $\in \mathbb{R}$,
\begin{align}
\label{def:T1}
\mathbb{P}_{0}\left(T<\infty\right)=\mathbb{P}_{1}\left(T<\infty\right)=1,
\end{align}
i.e., the random time $T$ defined in~\eqref{def:T} is a.s. finite under both the hypotheses. Indeed, we have,
\begin{align}
\label{eq:stop}
&\mathbb{P}_{1}\left(T > t\right)\le\mathbb{Q}\left(\frac{-\gamma_{d,i}^{h}+mt}{\sqrt{\frac{2mt}{N}+\frac{2mr^{2}(1-r^{2t})}{1-r^{2}}}}\right)\nonumber\\
&\Rightarrow\lim_{t\to\infty}\mathbb{P}_{1}\left(T > t\right)=0\nonumber\\
&\Rightarrow\mathbb{P}_{1}\left(T <\infty\right)=1.
\end{align}
The proof for $H_{0}$ follows in a similar way.

Now, since~\eqref{def:T1} holds, the quantity $S_{d,i}(T)$ is well-defined a.s. under $H_{0}$. Now, noting that, under $H_{0}$, for any $t$, the quantity $S_{d,i}(t)$ is Gaussian with mean $-mt$ and variance upper bounded by $\frac{2mt}{N}+\frac{2mr^{2}(1-r^{2t})}{1-r^{2}}$ (see Proposition~\ref{prop:SdGauss}), we have,
\begin{align}
\label{eq:3.59}
&\mathbb{P}_{FA}^{d,i}=\mathbb{P}_{0}(S_{d,i}(T)\ge \log\hat{B})=\sum_{t=1}^{\infty}\mathbb{P}_{0}(T=t,S_{d,i}(t)\ge \log\hat{B})\nonumber\\
&\le \sum_{t=1}^{\infty}\mathbb{P}_{0}(S_{d,i}(t)\ge \log\hat{B})\nonumber\\
&\leq\sum_{t=1}^{\infty}\mathbb{Q}\Big(\frac{\log\hat{B}+mt}{\sqrt{\frac{2mt}{N}+\frac{2mr^{2}(1-r^{2t})}{1-r^{2}}}}\Big).
\end{align}
To obtain a condition for $\gamma_{d,i}^{h}$ in the $\cisprt$ such that $\mathbb{P}_{FA}^{d,i}\leq\alpha$, let's define $k > 0 $ such that $k=Nr^{2}$. Now, note that $k$ thus defined satisfies
\begin{align}
\label{eq:3.6a}
&\frac{2mr^{2}(1-r^{2t})}{1-r^{2}}\le \frac{2mkt}{N}, ~~\forall t.
\end{align}

Then we have, by~\eqref{eq:3.59}-\eqref{eq:3.6a},
\allowdisplaybreaks[1]
\begin{align}
\label{eq:3.6}
&\mathbb{P}_{FA}^{d,i}\le \sum_{t=1}^{\infty}\mathbb{Q}\Big(\frac{\log\hat{B}+mt}{\sqrt{\frac{2mt}{N}+\frac{2mr^{2}(1-r^{2t})}{1-r^{2}}}}\Big)
\le \sum_{t=1}^{\infty}\mathbb{Q}\Big(\frac{\log\hat{B}+mt}{\sqrt{\frac{2mt(k+1)}{N}}}\Big)
\le \frac{1}{2}\sum_{t=1}^{\infty}e^{\frac{-(\gamma_{d,i}^{h})^{2}-m^{2}t^{2}-2\gamma_{d,i}^{h}mt}{\frac{4mt(k+1)}{N}}}\nonumber\\
&=\frac{e^{-\frac{N\gamma_{d,i}^{h}}{2(k+1)}}}{2}\Big(\sum_{t=1}^{\lfloor\frac{\gamma_{d,i}^{h}}{2m}\rfloor}e^{\frac{-N(\gamma_{d,i}^{h})^{2}-Nm^{2}t^{2}}{4mt(k+1)}}+\sum_{t=\lfloor\frac{\gamma_{d,i}^{h}}{2m}\rfloor+1}^{\lfloor\frac{\gamma_{d,i}^{h}}{m}\rfloor}e^{\frac{-N(\gamma_{d,i}^{h})^{2}-Nm^{2}t^{2}}{4mt(k+1)}}\nonumber\\
&+\sum_{t=\lfloor\frac{\gamma_{d,i}^{h}}{m}\rfloor+1}^{\lfloor\frac{2\gamma_{d,i}^{h}}{m}\rfloor}e^{\frac{-N(\gamma_{d,i}^{h})^{2}-Nm^{2}t^{2}}{4mt(k+1)}}+\sum_{t=\lfloor\frac{2\gamma_{d,i}^{h}}{m}\rfloor+1}^{\infty}e^{\frac{-N(\gamma_{d,i}^{h})^{2}-Nm^{2}t^{2}}{4mt(k+1)}}\Big)\nonumber\\
&\le \frac{e^{-\frac{N\gamma_{d,i}^{h}}{2(k+1)}}}{2}\Big( e^{-\frac{N\gamma_{d,i}^{h}}{2(k+1)}}\underbrace{\sum_{t=1}^{\lfloor\frac{\gamma_{d,i}^{h}}{2m}\rfloor}e^{\frac{-Nmt}{4(k+1)}}}_\text{(1)}+e^{-\frac{N\gamma_{d,i}^{h}}{4(k+1)}}\underbrace{\sum_{t=\lfloor\frac{\gamma_{d,i}^{h}}{2m}\rfloor+1}^{\lfloor\frac{\gamma_{d,i}^{h}}{m}\rfloor}e^{\frac{-Nmt}{4(k+1)}}}_\text{(2)}\nonumber\\
&+e^{-\frac{N\gamma_{d,i}^{h}}{8(k+1)}}\underbrace{\sum_{t=\lfloor\frac{\gamma_{d,i}^{h}}{m}\rfloor+1}^{\lfloor\frac{2\gamma_{d,i}^{h}}{m}\rfloor}e^{\frac{-Nmt}{4(k+1)}}}_\text{(3)}+\underbrace{\sum_{t=\lfloor\frac{2\gamma_{d,i}^{h}}{m}\rfloor+1}^{\infty}e^{\frac{-Nmt}{4(k+1)}}}_\text{(4)}\Big)\nonumber\\
&\le \frac{e^{-\frac{N\gamma_{d,i}^{h}}{2(k+1)}}}{2(1-e^{-\frac{Nm}{4(k+1)}})}\Big( e^{-\frac{N\gamma_{d,i}^{h}}{2(k+1)}}+e^{-\frac{N\gamma_{d,i}^{h}}{4(k+1)}}e^{-\frac{N\gamma_{d,i}^{h}}{8(k+1)}}+e^{-\frac{N\gamma_{d,i}^{h}}{8(k+1)}}e^{-\frac{N\gamma_{d,i}^{h}}{4(k+1)}}+e^{-\frac{N\gamma_{d,i}^{h}}{2(k+1)}}\Big)\nonumber\\
&\le  \frac{2e^{-\frac{7N\gamma_{d,i}^{h}}{8(k+1)}}}{1-e^{-\frac{Nm}{4(k+1)}}}.
\end{align}

\noindent In the above set of equations we use the fact that $\mathbb{Q}(x)$ is a non-increasing function, the inequality $\mathbb{Q}(x)\le\frac{1}{2}e^{\frac{-x^{2}}{2}}$, and we upper bound $(1)-(4)$ by their infinite geometric sums.

\noindent We now note that, a sufficient condition for $\mathbb{P}_{FA}^{d,i}\le\alpha$ to hold is the following:
\begin{align}
\label{eq:3.61}
\frac{2e^{-\frac{7N\gamma_{d,i}^{h}}{8(k+1)}}}{1-e^{-\frac{Nm}{4(k+1)}}}\le\alpha.
\end{align}
Solving \eqref{eq:3.61}, we have that, any $\gamma_{d,i}^{h}$ that satisfies
\begin{align}
\label{eq:3.62}
\gamma_{d,i}^{h}\ge \gamma_{d}^{h,0}=\frac{8(k+1)}{7N}\left(\log\left(\frac{2}{\alpha}\right)-\log(1-e^{-\frac{Nm}{4(k+1)}})\right),
\end{align}
achieves $\mathbb{P}_{FA}^{d,i}\le\alpha$ in the $\cisprt$.

Proceeding as in \eqref{eq:3.59} and \eqref{eq:3.6} we have that, any $\gamma_{d,i}^{l}$ that satisfies
\begin{align}
\label{eq:3.62a}
\gamma_{d,i}^{l}\le \gamma_{d}^{l,0}\doteq\frac{8(k+1)}{7N}\left(\log\left(\frac{\beta}{2}\right)+\log(1-e^{-\frac{Nm}{4(k+1)}})\right),
\end{align}
achieves $\mathbb{P}_{M}^{d,i}\le\beta$ in the $\cisprt$.

Clearly, by the above, any pair $(\gamma_{d,i}^{h},\gamma_{d,i}^{l})$ satisfying $\gamma_{d,i}^{h}\in [\gamma_{d}^{h,0},\infty)$ and $\gamma_{d,i}^{l}\in (-\infty,\gamma_{d}^{l,0}]$ (see~\eqref{eq:3.62} and~\eqref{eq:3.62a}) ensures that $\mathbb{P}_{FA}^{d,i}\le\alpha$ and $\mathbb{P}_{M}^{d,i}\leq\beta$. The a.s. finiteness of the corresponding stopping time $T_{d,i}$ (see~\eqref{eq:f3}) under both $H_{0}$ and $H_{1}$ follows readily by arguments as in~\eqref{def:T1}.
\end{IEEEproof}
\begin{Remark}
\label{rem:tighter}
It is to be noted that the derived thresholds are sufficient conditions only. The approximations (see $(1)-(4)$ in \eqref{eq:3.6}) made in the steps of deriving the expressions of the thresholds were done so as to get a tractable expression of the range. By solving the following set of equations
\begin{align}
\label{eq:3.12111}
&\frac{1}{2}\sum_{t=1}^{\infty}e^{\frac{-N(\gamma_{d,i}^{l})^{2}-Nm^{2}t^{2}+2N\gamma_{d,i}^{l}mt}{4mt(k+1)}}\le \beta \nonumber\\
&\frac{1}{2}\sum_{t=1}^{\infty}e^{\frac{-N(\gamma_{d,i}^{h})^{2}-Nm^{2}t^{2}-2N\gamma_{d,i}^{h}mt}{4mt(k+1)}}\le \alpha
\end{align}
numerically, tighter thresholds can be obtained.
\end{Remark}

\noindent\begin{IEEEproof}[Proof of Lemma \ref{1.2}]
Let us define the event $A_{s}^{i}$ as $\{\gamma_{d,i}^{l}\le S_{d,i}(s) \le \gamma_{d,i}^{h}\}$.
Now, note that
\begin{align}
\label{eq:12}
\mathbb{P}_{1}(T_{d,i}>t)=\mathbb{P}_{1}(\cap_{s=1}^{t}A_{s}^{i}),
\end{align}
and
\begin{align}
\label{eq:13}
\mathbb{P}_{1}(\cap_{s=1}^{t}A_{s}^{i}) \le \mathbb{P}_{1}(A_{t}^{i}).
\end{align}
By Proposition~\ref{prop:SdGauss}, under $H_{1}$, for any $t$, the quantity $S_{d,i}(t)$ is Gaussian with mean $mt$ and variance upper bounded by $\frac{2mt}{N}+\frac{2mr^{2}(1-r^{2t})}{1-r^{2}}$. Hence we have, for all $i=1,2,\hdots,N$.
\begin{align}
\label{eq:14}
&\mathbb{P}_{1}(T_{d,i}>t) \le\mathbb{Q}\Big(\frac{-\gamma_{d,i}^{h}+mt}{\sqrt{\frac{2mt}{N}+\frac{2mr^{2}(1-r^{2t})}{1-r^{2}}}}\Big).
\end{align}
\end{IEEEproof}
\noindent\begin{IEEEproof}[Proof of Corollary \ref{7.2}]
For simplicity of notation, let $a=\frac{Nm}{4}$ and $b=\frac{\sigma^{2}\pi^{2}}{2N(\gamma_{c}^{h}-\gamma_{c}^{l})^{2}}$.
From \eqref{eq:1.01}, we have,
\begin{align}
\label{eq:12.6}
&\frac{1}{t}\log(\mathbb{P}_{1}(T_{c}> t))\geq\frac{1}{t}\log\left(\exp\left(\frac{N\mu\gamma_{c}^{l}}{\sigma^{2}}\right)K_{t}^{\infty}\left(\gamma_{c}^{h}\right)-\exp\left(\frac{N\mu\gamma_{c}^{h}}{\sigma^{2}}\right)K_{t}^{\infty}\left(\gamma_{c}^{l}\right)\right)\nonumber\\
&=\frac{1}{t}\log\left(\exp\left(-\left(a+b\right)t\right)\right)\nonumber\\&+\frac{1}{t}\log\left(b\sum_{s=1}^{\infty}\frac{s(-1)^{s+1}}{a+s^{2}b}\exp\left(-b(s^{2}-1)t\right)\right.\nonumber\\&\left.\times\left(\exp\left(\frac{N\mu\gamma_{c}^{l}}{\sigma^{2}}\right)\sin\left(\frac{s\pi\gamma_{c}^{h}}{\gamma_{c}^{h}-\gamma_{c}^{l}}\right)-\exp\left(\frac{N\mu\gamma_{c}^{h}}{\sigma^{2}}\right)\sin\left(\frac{s\pi\gamma_{c}^{l}}{\gamma_{c}^{h}-\gamma_{c}^{l}}\right)\right)\right).
\end{align}

For all $t,S\geq 1$, let
\begin{align}
U(t,S)=\frac{1}{t}\log\left(b\sum_{s=1}^{S}\frac{s(-1)^{s+1}}{a+s^{2}b}\exp\left(-b(s^{2}-1)t\right)\right.\nonumber\\\left.\times\left(\exp\left(\frac{N\mu\gamma_{c}^{l}}{\sigma^{2}}\right)\sin\left(\frac{s\pi\gamma_{c}^{h}}{\gamma_{c}^{h}-\gamma_{c}^{l}}\right)-\exp\left(\frac{N\mu\gamma_{c}^{h}}{\sigma^{2}}\right)\sin\left(\frac{s\pi\gamma_{c}^{l}}{\gamma_{c}^{h}-\gamma_{c}^{l}}\right)\right)\right)
\end{align}
and let $g=\exp\left(\frac{N\mu\gamma_{c}^{h}}{\sigma^{2}}\right)+\exp\left(\frac{N\mu\gamma_{c}^{l}}{\sigma^{2}}\right)$.

Note that for all $t\geq 1$, the limit
\begin{equation}
\label{eq:12345}
\lim_{S\rightarrow\infty}U(t,S)
\end{equation}
 exists and is finite (by Theorem~\ref{1.1}), and similarly for all $S\geq 1$,
\begin{align}
\label{eq:1234}
\lim_{t\rightarrow\infty}U(t,S)=\lim_{t\to\infty}\frac{1}{t}\log\left(b\sum_{s=1}^{S}\frac{s(-1)^{s+1}}{a+s^{2}b}\exp\left(-b(s^{2}-1)t\right)\right.\nonumber\\\left.\times\sin\left(\frac{s\pi\gamma_{c}^{h}}{\gamma_{c}^{h}-\gamma_{c}^{l}}\right)\left(\exp\left(\frac{N\mu\gamma_{c}^{l}}{\sigma^{2}}\right)+(-1)^{s+1}\exp\left(\frac{N\mu\gamma_{c}^{h}}{\sigma^{2}}\right)\right)\right)\nonumber\\
=\lim_{t\to\infty}\frac{1}{t}\log\left(\frac{bg}{a+b}\sin\left(\frac{\pi\gamma_{c}^{h}}{\gamma_{c}^{h}-\gamma_{c}^{l}}\right)\right)\nonumber\\ = 0,
\end{align}
where we use the fact that only the largest exponent in a finite summation of exponential terms contributes to its log-normalised limit as $t\rightarrow\infty$ and
\begin{equation}
\sin\left(\frac{s\pi\gamma_{c}^{h}}{\gamma_{c}^{h}-\gamma_{c}^{l}}\right)=(-1)^{s}\sin\left(\frac{s\pi\gamma_{c}^{l}}{\gamma_{c}^{h}-\gamma_{c}^{l}}\right).
\end{equation}
Finally, using the fact that there exists a constant $c_{5}>0$ (independent of $t$ and $S$) such that for all $t,S\geq 1$,
\begin{align}
\label{eq:12.61}
U(t,S)\le \frac{1}{t}\log\left(bg\sum_{s=1}^{S}\frac{s}{a+s^{2}b}\exp\left(-b(s^{2}-1)t\right)\right)\leq c_{5},
\end{align}
we may conclude that the convergence in~\eqref{eq:12345}-\eqref{eq:1234} are uniform in $S$ and $t$ respectively. This in turn implies that the order of the limits may be interchanged and we have that
\begin{equation}
\label{eq:123456}
\lim_{t\rightarrow\infty}\lim_{S\rightarrow\infty}U(t,S)=\lim_{S\rightarrow\infty}\lim_{t\rightarrow\infty}U(t,S)=0.
\end{equation}

Hence, we have from \eqref{eq:12.6} and~\eqref{eq:123456},
\begin{align}
\label{eq:12.62}
&\liminf_{t\to\infty}\frac{1}{t}\log(\mathbb{P}_{1}(T_{c}> t))\geq -(a+b)\nonumber\\&+\lim_{t\to\infty}\lim_{S\to\infty}\frac{1}{t}\log\left(b\sum_{s=1}^{S}\frac{s(-1)^{s+1}}{a+s^{2}b}\exp\left(-b(s^{2}-1)t\right)\right.\nonumber\\&\left.\times\sin\left(\frac{s\pi\gamma_{c}^{h}}{\gamma_{c}^{h}-\gamma_{c}^{l}}\right)\left(\exp\left(\frac{N\mu\gamma_{c}^{l}}{\sigma^{2}}\right)-(-1)^{s}\exp\left(\frac{N\mu\gamma_{c}^{h}}{\sigma^{2}}\right)\right)\right)\nonumber\\
&=-(a+b)+\lim_{S\to\infty}\lim_{t\to\infty}\frac{1}{t}\log\left(b\sum_{s=1}^{S}\frac{s(-1)^{s+1}}{a+s^{2}b}\exp\left(-b(s^{2}-1)t\right)\right.\nonumber\\&\left.\times\sin\left(\frac{s\pi\gamma_{c}^{h}}{\gamma_{c}^{h}-\gamma_{c}^{l}}\right)\left(\exp\left(\frac{N\mu\gamma_{c}^{l}}{\sigma^{2}}\right)+(-1)^{s+1}\exp\left(\frac{N\mu\gamma_{c}^{h}}{\sigma^{2}}\right)\right)\right)\nonumber\\
&=-(a+b)+\lim_{S\to\infty}\lim_{t\to\infty}\frac{1}{t}\log\left(\frac{bg}{a+b}\sin\left(\frac{\pi\gamma_{c}^{h}}{\gamma_{c}^{h}-\gamma_{c}^{l}}\right)\right)\nonumber\\
&=-(a+b)=-\left(\frac{Nm}{4}+\frac{\sigma^{2}\pi^{2}}{2N(\gamma_{c}^{h}-\gamma_{c}^{l})^{2}}\right).
\end{align}

\end{IEEEproof}
\noindent\begin{IEEEproof}[Proof of Theorem \ref{7.3}]
	We use the following upper bound for $\mathbb{Q}$ function in the proof below
	\begin{align}
		\label{eq:12.79}
		\mathbb{Q}(x)\le \frac{1}{x\sqrt{2\pi}}e^{-x^{2}/2}
	\end{align}
	From \eqref{eq:1.02},\eqref{eq:12.79} and \eqref{eq:12.61}, we have,
	{\small
		\begin{align}
			\label{eq:12.8}
			&\limsup_{t\to\infty}\frac{1}{t}\log(\mathbb{P}_{1}(T_{d,i}> t))\nonumber\\& \le \limsup_{t\to\infty}\frac{1}{t}\log\Big(\frac{1}{\sqrt{2\pi}}\frac{\sqrt{\frac{2mt}{N}+2m\frac{r^{2}(1-r^{2t})}{(1-r^{2})}}}{(-\gamma_{d,i}^{h}+mt)}e^{\frac{-N(-\gamma_{d,i}^{h}+mt)^{2}}{4mt+4mN\frac{r^{2}(1-r^{2t})}{(1-r^{2})}}}\Big)\nonumber\\
			&\le \limsup_{t\to\infty}\frac{1}{t}\left(\log\Big(\frac{\sqrt{\frac{2mt}{N}+2m\frac{r^{2}(1-r^{2t})}{(1-r^{2})}}}{\sqrt{2\pi}(mt-\gamma_{d,i}^{h})}\Big)-\frac{N(\gamma_{d,i}^{h})^{2}}{4mt+4m\frac{r^{2}(1-r^{2t})}{(1-r^{2})}}\right.\nonumber\\& \left.-\frac{Nmt}{4+4\frac{r^{2}(1-r^{2t})}{(t(1-r^{2}))}}+\frac{Nm\gamma_{d,i}^{h}t}{2mt+2mN\frac{r^{2}(1-r^{2t})}{(1-r^{2})}}\right)\nonumber\\
			&\Rightarrow\limsup_{t\to\infty}\frac{1}{t}\log(\mathbb{P}_{1}(T_{d,i}> t))\le -\frac{Nm}{4}.\nonumber\\
		\end{align}
	}
\end{IEEEproof}

\noindent The proof of Theorem~\ref{edct_1} requires an intermediate result that estimates the divergence between the agent statistics over time.
\noindent\begin{Lemma}
\label{lm:div_est}
Let the Assumptions~\ref{as:1}, \ref{as:3} and \ref{as:4} hold. Then, there exists a constant $c_{1}$, depending on the network topology and the Gaussian model statistics only, such that
\begin{equation}
\label{lm:div_est1}
\mathbb{E}_{1}\left[\sup_{t\geq 0}\|S_{d,i}(t)-S_{d,j}(t)\|\right]\leq c_{1}
\end{equation}
for all agent pairs $(i,j)$.
\end{Lemma}
\noindent\begin{IEEEproof}
Denoting by $\mathbf{S}_{d}(t) = t\mathbf{P}_{d}(t)$ the vector of the agent test statistics $S_{d,i}(t)$'s, we have by~\eqref{eq:3},
\begin{equation}
\label{lm:div_est2}
\mathbf{S}_{d}(t+1) = W\left(\mathbf{S}_{d}(t) + \mathbf{\eta}(t+1)\right).
\end{equation}
Let $\overline{S}_{d}(t)$ denote the average of the $S_{d,i}(t)$'s, i.e.,
\begin{equation}
\label{lm:div_est3}
\overline{S}_{d}(t)=\left(1/N\right).\left(S_{d,1}(t)+\cdots+ S_{d,N}(t)\right),
\end{equation}
Noting that $J\mathbf{S}_{d}(t)=\overline{S}_{d}(t)\mathbf{1}$ and $WJ=JW=J$, we have from~\eqref{lm:div_est2}
\begin{equation}
\label{lm:div_est4}
\mathbf{v}_{t+1} = \left(W-J\right)\mathbf{v}_{t}+\mathbf{u}_{t+1},
\end{equation}
where $\mathbf{v}_{t}$ and $\mathbf{u}_{t}$, for all $t\geq 0$, are given by
\begin{equation}
\label{lm:div_est5}
\mathbf{v}_{t} = \mathbf{S}_{d}(t) - \overline{S}_{d}(t)\mathbf{1}
\end{equation}
and
\begin{equation}
\label{lm:div_est6}
\mathbf{u}_{t+1} = \left(W-J\right)\mathbf{\eta}(t+1).
\end{equation}
It is important to note that the sequence $\{\mathbf{u}_{t}\}$ is i.i.d. Gaussian and, in particular, there exists a constant $c_{2}$ such that $\mathbb{E}_{1}[\|\mathbf{u}_{t}\|^{2}]\leq c_{2}$ for all $t$.

Now, by~\eqref{lm:div_est4} we obtain
\begin{equation}
\label{lm:div_est7}
\|\mathbf{v}_{t+1}\|\leq r\|\mathbf{v}_{t}\|+\|\mathbf{u}_{t+1}\|,
\end{equation}
where recall $r=\|W-J\|<1$. Since the sequence $\{\mathbf{u}_{t}\}$ is i.i.d. and $\mathcal{L}_{2}$-bounded, an application of the Robbins-Siegmund's lemma (see~\cite{baldi2002martingales}) yields
\begin{equation}
\label{lm:div_est8}
\mathbb{E}_{1}\left[\sup_{t\geq 0}\|\mathbf{v}_{t}\|\right]\leq c_{3}<\infty,
\end{equation}
where $c_{3}$ is a constant that may be chosen as a function of $r$, $c_{2}$ and $\mathbb{E}_{1}[\|\mathbf{v}_{0}\|]$. Now, noting that, for any pair $(i,j)$,
\begin{align}
\label{lm:div_est9}
\mathbb{E}_{1}\left[\sup_{t\geq 0}\|S_{d,i}(t)-S_{d,j}(t)\|\right]\leq \mathbb{E}_{1}\left[\sup_{t\geq 0}\|S_{d,i}(t)-\overline{S}_{d}(t)\|\right]+\mathbb{E}_{1}\left[\sup_{t\geq 0}\|S_{d,j}(t)-\overline{S}_{d}(t)\|\right]\leq 2c_{3},
\end{align}
the desired assertion follows.
\end{IEEEproof}

\noindent\begin{IEEEproof}[Proof of Theorem \ref{edct_1}] We prove the upper bound in Theorem~\ref{edct_1} first. Since $\mathbb{P}_{1}(T_{d,i}<\infty)=1$, for the upper bound we have,
\allowdisplaybreaks[1]
\begin{align}
\label{eq:edcp211}
&\mathbb{E}_{1}[T_{d,i}]=\sum_{t=0}^{\infty}\mathbb{P}_{1}(T_{d,i} > t)\nonumber\\&\overset{(a)}{\le}\sum_{0}^{\infty}\mathbb{Q}\big(\frac{-\gamma_{d,i}^{h}+mt}{\sqrt{\frac{2mt}{N}+\frac{2mr^{2}(1-r^{2t})}{1-r^{2}}}}\Big)\nonumber\\
&=\underbrace{\sum_{0}^{\lfloor\frac{\gamma_{d,i}^{h}}{m}\rfloor}\mathbb{Q}\Big(\frac{-\gamma_{d,i}^{h}+mt}{\sqrt{\frac{2mt}{N}+\frac{2mr^{2}(1-r^{2t})}{1-r^{2}}}}\Big)}_\text{(1)}+\underbrace{\sum_{\lfloor\frac{\gamma_{d,i}^{h}}{m}\rfloor+1}^{\lfloor\frac{3\gamma_{d,i}^{h}}{2m}\rfloor}\mathbb{Q}\Big(\frac{-\gamma_{d,i}^{h}+mt}{\sqrt{\frac{2mt}{N}+\frac{2mr^{2}(1-r^{2t})}{1-r^{2}}}}\Big)}_\text{(2)}\nonumber\\&+\underbrace{\sum_{\lfloor\frac{3\gamma_{d,i}^{h}}{2m}\rfloor+1}^{\lfloor\frac{2\gamma_{d,i}^{h}}{m}\rfloor}\mathbb{Q}\Big(\frac{-\gamma_{d,i}^{h}+mt}{\sqrt{\frac{2mt}{N}+\frac{2mr^{2}(1-r^{2t})}{1-r^{2}}}}\Big)}_\text{(3)}+\underbrace{\sum_{\lfloor\frac{2\gamma_{d,i}^{h}}{m}\rfloor+1}^{\infty}\mathbb{Q}\Big(\frac{-\gamma_{d,i}^{h}+mt}{\sqrt{\frac{2mt}{N}+\frac{2mr^{2}(1-r^{2t})}{1-r^{2}}}}\Big)}_\text{(4)}\nonumber\\
&\overset{(b)}{\le} \frac{\gamma_{d,i}^{h}}{m}+\frac{\gamma_{d,i}^{h}}{4m}+\frac{1}{2}e^{\frac{N\gamma_{d,i}^{h}}{2(k+1)}}\sum_{\lfloor\frac{3\gamma_{d,i}^{h}}{2m}\rfloor+1}^{\lfloor\frac{2\gamma_{d,i}^{h}}{m}\rfloor}e^{\frac{-(N\gamma_{d,i}^{h})^{2}-Nm^{2}t^{2}}{4m(k+1)t}}+\frac{1}{2(1-e^{\frac{-Nm}{4(k+1)}})}\nonumber\\
&\le \frac{5\gamma_{d,i}^{h}}{4m}+\frac{1}{2(1-e^{\frac{-Nm}{4(k+1)}})}+\frac{1}{2}e^{\frac{3N\gamma_{d,i}^{h}}{8(k+1)}}\sum_{\lfloor\frac{3\gamma_{d,i}^{h}}{2m}\rfloor+1}^{\lfloor\frac{2\gamma_{d,i}^{h}}{m}\rfloor}e^{\frac{-Nmt}{4(k+1)}}\nonumber\\
&\le \frac{5\gamma_{d,i}^{h}}{4m}+\frac{1}{1-e^{\frac{-Nm}{4(k+1)}}},
\end{align}
where $(a)$ is due to the upper bound derived in Lemma \ref{1.2} and $(b)$ is due to the following : 1) $\forall t \in [0,\lfloor\frac{\gamma_{d,i}^{h}}{m}\rfloor]$ in $(1)$, $-\gamma_{d,i}^{h}+mt$ is negative and hence every term in the summation can be upper bounded by $1$; 2) $\forall t \in [\lfloor\frac{\gamma_{d,i}^{h}}{m}\rfloor+1, \lfloor\frac{3\gamma_{d,i}^{h}}{2m}\rfloor ]$ in $(2)$, $-\gamma_{d,i}^{h}+mt$ is positive and hence every term in the summation can be upper bounded by $\frac{1}{2}$; and 3) for the terms $(3)$ and $(4)$, the inequality $\mathbb{Q}(x) \le \frac{1}{2}e^{-x^{2}/2}$ is used and the sums are upper bounded by summing the resulting geometric series.

In order to obtain the lower bound, we first note that conditioned on hypothesis $H_{1}$, at the stopping time $T_{d,i}$,  an agent exceeds the threshold  $\gamma_{d,i}^{h}$ with probability at least $1-\epsilon$ and is lower than the threshold $\gamma_{d,i}^{l}$ with probability at most $\epsilon$. Moreover, with $\alpha=\beta=\epsilon$, $\gamma_{d,i}^{h}=-\gamma_{d,i}^{l}$.

Now, denote by $E^{h}_{i}$ the event $E^{h}_{i} = \{S_{d,i}(T_{d,i}) \geq \gamma_{d,i}^{h}\}$ and by $E^{l}_{i}$ the event $E^{l}_{i} = \{S_{d,i}(T_{d,i}) \leq \gamma_{d,i}^{l}\}$. Since $\mathbb{P}_{1}(T_{d,i} < \infty) =1$, we have that
\begin{equation}
\label{ext1}
\mathbb{E}_{1}\left[S_{d,i}(t)\right]=\mathbb{E}_{1}\left[S_{d,i}(t).\mathbb{I}_{E^{h}_{i}}\right]+ \mathbb{E}_{1}\left[S_{d,i}(t).\mathbb{I}_{E^{l}_{i}}\right],
\end{equation}
where $\mathbb{I}_{\{\cdot\}}$ denotes the indicator function. We now lower bound the quantities on the R.H.S. of~\eqref{ext1}. Note that $\gamma_{d,i}^{h}\geq 0$ and  $S_{d,i}(t)\geq\gamma_{d,i}^{h}$ on $E^{h}_{i}$. Hence
\begin{equation}
\label{ext2}
\mathbb{E}_{1}\left[S_{d,i}(t).\mathbb{I}_{E^{h}_{i}}\right]\geq \gamma_{d,i}^{h}\mathbb{P}_{1}\left(E^{h}_{i}\right)\geq (1-\epsilon)\gamma_{d,i}^{h}.
\end{equation}
Now recall the construction in the proof of Lemma~\ref{lm:div_est} and note that by~\eqref{lm:div_est2} we have
\begin{equation}
\label{ext3}
S_{d,i}(t) = S_{d,i}(t-1)-\sum_{j\in\Omega_{i}}w_{ij}\left(S_{d,i}(t-1)-S_{d,j}(t-1)\right)+\eta_{i}(t).
\end{equation}
Hence, we have that
\begin{align}
\label{ext40}
S_{d,i}(T_{d,i}).\mathbb{I}_{E^{l}_{i}} \\ \geq S_{d,i}(T_{d,i}-1).\mathbb{I}_{E^{l}_{i}} - \sum_{j\in\Omega_{i}}w_{ij}\|S_{d,i}(T_{d,i}-1)-S_{d,j}(T_{d,i}-1)\|-\|\eta_{i}(T_{d,i})\|\\
\geq S_{d,i}(T_{d,i}-1).\mathbb{I}_{E^{l}_{i}} - \sum_{j\in\Omega_{i}}w_{ij}\sup_{t\geq 0}\|S_{d,i}(t)-S_{d,j}(t)\|-\|\eta_{i}(T_{d,i})\|.
\end{align}
Now, observe that on the event $E^{l}_{i}$, $S_{d,i}(T_{d,i}-1)> \gamma_{d,i}^{l}$ a.s. Since $\gamma_{d,i}^{l}<0$ and $\mathbb{P}_{1}(E^{l}_{i})\leq\epsilon$ (by hypothesis), we have that
\begin{align}
\label{ext4}
\gamma_{d,i}^{l}\epsilon\leq\gamma_{d,i}^{l}\mathbb{P}_{1}(E^{l}_{i}) \\ = \mathbb{E}_{1}\left[\gamma_{d,i}^{l}.\mathbb{I}_{E^{l}_{i}}\right]\leq \mathbb{E}_{1}\left[S_{d,i}(T_{d,i}-1).\mathbb{I}_{E^{l}_{i}}\right].
\end{align}
Note that, by Lemma~\ref{lm:div_est}, we have
\begin{align}
\label{ext5}
\mathbb{E}_{1}\left[\sum_{j\in\Omega_{i}}w_{ij}\sup_{t\geq 0}\|S_{d,i}(t)-S_{d,j}(t)\|\right]\\
\leq \sum_{j\in\Omega_{i}}w_{ij}\mathbb{E}_{1}\left[\sup_{t\geq 0}\|S_{d,i}(t)-S_{d,j}(t)\|\right]\\
\leq |\Omega_{i}|c_{1}.
\end{align}
Finally, by arguments similar to~\cite{wald1973sequential,lorden1970excess} for characterizing expected overshoots in stopped random sums (see, in particular, Theorem 1 in~\cite{lorden1970excess}) it follows that there exists a constant $c_{4}$ (depending on the Gaussian model statistics and the network topology only) such that
\begin{equation}
\label{ext6}
\mathbb{E}_{1}\left[\|\eta_{i}(T_{d,i})\|\right]\leq c_{4}.
\end{equation}
In particular, note that, the constant $c_{4}$ in~\eqref{ext6} may be chosen to be independent of the thresholds and, hence, the error tolerance parameter $\epsilon$.
Substituting~\eqref{ext4}-\eqref{ext6} in~\eqref{ext40} we obtain
\begin{align}
\label{ext7}
\mathbb{E}_{1}\left[S_{d,i}(T_{d,i}).\mathbb{I}_{E^{l}_{i}}\right]\geq \gamma_{d,i}^{l}\epsilon - |\Omega_{i}|c_{1} - c_{4}.
\end{align}
This together with~\eqref{ext1}-\eqref{ext2} yield
\begin{align}
\label{ext8}
\mathbb{E}_{1}\left[S_{d,i}(T_{d,i})\right]\geq\left(1-\epsilon\right)\gamma_{d,i}^{h}+\gamma_{d,i}^{l}\epsilon - |\Omega_{i}|c_{1} - c_{4}\\
=\left(1-2\epsilon\right)\gamma_{d,i}^{h} - c,
\end{align}
where the last equality follows by noting that $\gamma_{d,i}^{h}=-\gamma_{d,i}^{l}$ and taking the constant $c$ to $c=|\Omega_{i}|c_{1} + c_{4}$.

We note that the event $\{T_{d,i}=t\}$ is independent of $\eta_{i}, i > t$. We also have from Theorem \ref{dist_th} that $\mathbb{P}_{1}(T_{d,i} < \infty) =1$. Hence, we have,
\begin{align}
\label{eq:edcp11a}
&\mathbb{E}_{1}[S_{d,i}(T_{d,i})]=\mathbb{E}_{1}[\sum_{j=1}^{T_{d,i}}\mathbf{e_{i}}^{\top}\mathbf{W}^{t+1-j}\mathbf{\eta}(j)]\nonumber\\
&=\mathbb{E}_{1}\left[\sum_{j=1}^{\infty}\mathbb{I}_{\left\{T_{d,i}\ge j\right\}}\mathbf{e_{i}}^{\top}\mathbf{W}^{T_{d,i}+1-j}\mathbf{\eta}(j)\right]\nonumber\\
&=\sum_{j=1}^{\infty}\mathbb{E}_{1}\left[\mathbb{I}_{\left\{T_{d,i}\ge j\right\}}\mathbf{e_{i}}^{\top}\mathbf{W}^{T_{d,i}+1-j}\right]\mathbb{E}_{1}\left[\mathbf{\eta}(j)\right]\nonumber\\
&=m\sum_{j=1}^{\infty}\mathbb{E}_{1}\left[\mathbb{I}_{\left\{T_{d,i}\ge j\right\}}\mathbf{e_{i}}^{\top}\mathbf{W}^{T_{d,i}+1-j}\right]\mathbf{1}\nonumber\\
&=m\sum_{j=1}^{\infty}\mathbb{E}_{1}\left[\mathbb{I}_{\left\{T_{d,i}\ge j\right\}}\mathbf{e_{i}}^{\top}\mathbf{W}^{T_{d,i}+1-j}\mathbf{1}\right]\nonumber\\
&=m\mathbb{E}_{1}\left[T_{d,i}\right].
\end{align}
Combining \eqref{eq:edcp11a} and \eqref{ext8} we have,
\begin{align}
\label{eq:edcp11b}
\frac{(1-2\epsilon)\gamma_{d,i}^{h}}{m} - \frac{c}{m}\le \mathbb{E}_{1}[T_{d,i}]
\end{align}
and the desired assertion follows.

\end{IEEEproof}

\begin{IEEEproof}[Proof of Theorem \ref{edct}]
From \eqref{eq:edc1}, we first note that
\begin{align}
\label{eq:edcp1}
\frac{\mathbb{E}_{1}[T_{d,i}]}{\mathbb{E}_{1}[T_{c}]}\ge 1,\forall i=1,2,\hdots,N.
\end{align}
From the upper bound for the stopping time distribution derived for the  $\cisprt$  in \eqref{eq:14}, we have the following upper bound for $\mathbb{E}_{1}[T_{d,i}]$
\begin{align}
\label{eq:edcp2}
\mathbb{E}_{1}[T_{d,i}]\le \frac{5\gamma_{d,i}^{h}}{4m}+\frac{1}{1-e^{\frac{-Nm}{4(k+1)}}}.
\end{align}
We choose the threshold $\gamma_{d,i}^{h}$ to be
\begin{align}
\label{eq:edcp2b}
\gamma_{d,i}^{h}=\gamma_{d}^{h,0}=\frac{8(k+1)}{7N}(\log(\frac{2}{\epsilon})-\log(1-e^{\frac{-Nm}{4(k+1)}})).
\end{align}
Using \eqref{eq:edcp2} and \eqref{eq:edcp2b}, we have
\begin{align}
\label{eq:edcp3}
\limsup_{\epsilon\to0}\frac{\mathbb{E}_{1}[T_{d,i}]}{\mathbb{E}_{1}[T_{c}]}\le\lim_{\epsilon\to0}\frac{\frac{10}{7}(k+1) \log(\frac{2}{\epsilon})+O(1)}{(1-2\epsilon)\log(\frac{1-\epsilon}{\epsilon})}.
\end{align}
Noting that,
\begin{align}
\label{eq:edcp4}
\limsup_{\epsilon\to0}\frac{O(1)}{(1-2\epsilon)\log(\frac{1-\epsilon}{\epsilon})}= 0,
\end{align}
we obtain
\begin{align}
\label{eq:edcp5}
\limsup_{\epsilon\to0}\frac{\mathbb{E}_{1}[T_{d,i}]}{\mathbb{E}_{1}[T_{c}]}\le \frac{10(k+1)}{7}.
\end{align}
Combining \eqref{eq:edcp5} and \eqref{eq:edcp1}, the result follows.
\end{IEEEproof}
\section{Conclusion}
\label{sec:conc}
In this paper we have considered sequential detection of Gaussian binary hypothesis observed by a sparsely interconnected network of agents. The $\cisprt$ algorithm we proposed combines two terms : a \emph{consensus} term that updates at
each sensor its test statistic with the test statistics provided by agents in its one-hop neighborhood and an \emph{innovation} term that updates the current agent test statistic with the new local sensed information. We have shown that the $\cisprt$ can be designed to achieve a.s. finite stopping at each network agent with guaranteed error performance. We have provided explicit characterization of the large deviation decay exponents of tail probabilities of the $\cisprt$ stopping and its expected stopping time as a function of the network connectivity. The performance of the $\cisprt$ was further benchmarked w.r.t. the optimal centralized sequential detector, the SPRT. The techniques developed in this paper are of independent interest and we envision their applicability to other distributed sequential procedures. An interesting future direction would be to consider networks with random time-varying topology. We also intend to develop extensions of the $\cisprt$ for setups with correlated and non-linear non-Gaussian observation models.

\section{Acknowledgment}
\label{sec:ack}
The authors would like to thank the associate editor and the reviewers for their comments and detailed feedback that helped to improve the clarity and content of the paper and, in particular, an anonymous reviewer for pointing out a technical gap in a previous version of Corollary 4.5.
\bibliographystyle{IEEEtran}
\bibliography{dsprt,CentralBib}
\end{document}